\documentclass[10pt,twoside,reqno]{article}
\usepackage{amssymb,amsmath,latexsym}
\usepackage[varg]{pxfonts}
\newtheorem{thm}{Theorem}[section]

\newtheorem{lem}[thm]{Lemma}
\newtheorem{prop}[thm]{Proposition}
\newtheorem{defn}[thm]{Definition}
\newtheorem{remark}[thm]{Remark}
\newtheorem{alg}[thm]{Algorithm}
\newtheorem{example}[thm]{Example}
\numberwithin{equation}{section}

\begin{document}
\setcounter{page}{1}
\begin{center}
\vspace{0.4cm} {\large{\bf Generalized variational inclusion governed by generalized $\alpha\beta$-$H((., .), (., .))$-mixed accretive mapping in real\\ $q$-uniformly smooth Banach spaces}}\\
\vspace{0.4cm}
Sanjeev Gupta$^{a}$, Shamshad Husain$^b$ Vishnu Narayan Mishra$^{c\dag}$\\~~\\
$^a$Department of Humanities and Social Sciences\\
Indian Institute of Technology Kanpur, Kanpur-208016, India\\~\\
$^b$Department of Applied Mathematics\\
Faculty of Engineering ${\&}$ Technology\\
Aligarh Muslim University, Aligarh-202002, India\\
~~~\\
$^{c}$Department of Applied Mathematics ${\&}$ Humanities\\
Sardar Vallabhbhai National Institute of Technology, \\
Surat-395 007, India\\
\let\thefootnote\relax\footnotetext{$^c{\dag}$Correspondence}
\let\thefootnote\relax\footnotetext{Email address: s$_{-}$husain68@yahoo.com (S. Husain),~guptasanmp@gmail.com (S. Gupta),}
\let\thefootnote\relax\footnotetext{\hspace{2cm} vishnunarayanmishra@gmail.com (V.N. Mishra)}
\end{center}

\begin{abstract} In this paper, we investigate a new notion of accretive mappings called
generalized $\alpha\beta$-$H((.,.),(.,.))$-mixed accretive mappings in Banach spaces. We extend the concept of proximal-point mappings
associated with generalized $m$-accretive mappings to the
generalized $\alpha\beta$-$H((.,.),(.,.))$-mixed accretive mappings and prove that the
proximal-point mapping associated with generalized $\alpha\beta$-$H((.,.),(.,.))$-mixed accretive mapping is
single-valued and Lipschitz continuous. Some examples are given to
justify the definition of generalized $\alpha\beta$-$H((.,.),(.,.))$-mixed accretive
mappings. Further, by using the proximal mapping technique, an iterative algorithm for solving
a class of variational inclusions is constructed in real $q$-uniformly smooth Banach spaces. Under some suitable
conditions, we prove the convergence of iterative sequence generated by the algorithm.\\


\end{abstract}
\section{Introduction}

\par
\indent~~~~Variational inclusions, as the generalization of variational inequalities, have been widely studied in recent years. One of
the most interesting and important problems in the theory of variational inclusions is the development of an efficient and
implementable iterative algorithm. Variational inclusions include variational, quasi-variational, variational-like inequalities
as special cases. For application of variational inclusions, see for example \cite{AU}. Various kinds of iterative methods have been
studied to find the approximate solutions for variational inclusions. Among these methods, The proximal-point mapping techniques are important to study the
existence of solutions and to design iterative schemes for different
kinds of variational inequalities and their generalizations, which
are providing mathematical models to some problems arising in
optimization and control, economics, and engineering sciences, has been widely used by many authors. For details, we refer to see \cite{A-D}-\cite{Luo}, \cite{Peng1}-\cite{X-W},\cite{Z-H,Z-H1} and the references therein.\\
%
\par In order to study various variational inequalities and variational
inclusions, Huang and Fang \cite{H-F1} were the first to introduce the generalized $m$-accretive mapping
and give the definition of the proximal-point mapping in Banach spaces. Since then a number of researchers
introduced several classes of generalized $m$-accretive mappings such as $H$-accretive, $(H,~\eta)$-accretive, $(A,~\eta)$-accretive, $(P,~\eta)$-accretive mappings, see for examples \cite{F-H1}-\cite{F-H-T}, \cite{K-F1,K-F2,L-C-V,Peng1,Peng2}. Sun et al. \cite{Sun} introduced a new class of $M$-monotone mapping in Hilbert spaces. Recently, Zou and Huang \cite{Z-H,Z-H1}, Kazmi et al. \cite{K-A-S,K-B-A} introduced and studied a class of $H(.,.)$-accretive mappings, Ahmad and Dilshad \cite{A-D} introduced and studied a class of $H(.,.)$-$\eta$-cocoercive and Husain and Gupta \cite{H-G2} introduced and studied a class of $H((.,.),(.,.))$-mixed cocoercive mappings in Banach (Hilbert) spaces, an natural extension of $M$-monotone mapping and studied variational inclusions involving these mappings. In recent past, the methods based on different classes of proximal-point mappings have been developed to study
the existence of solutions and to discuss convergence and stability analysis of iterative algorithms for various classes of variational
inclusions, see for example \cite{A-D,B-H-F}, \cite{H-F1}-\cite{Luo}, \cite{Peng1}-\cite{X-W}, \cite{Z-H}.
\\

\par Very recently, Luo and Huang \cite{Luo} introduced and studied a class of $B$-monotone and Kazmi et al. \cite{K-A-S} introduced and studied a class of generalized $H(.,.)$-accretive mappings in Banach spaces, an extension of $H$-monotone mappings \cite{F-H1}. They showed some properties of the proximal-point mapping associated with $B$-monotone and generalized $H(.,.)$-accretive mapping, and obtained some applications for solving variational inclusions in Banach spaces.\\

\par Motivated and inspired by the research works mentioned above, we consider a new
class of generalized $\alpha\beta$-$H((.,.),(.,.))$-mixed accretive mappings for solving generalized set-valued variational inclusions in real $q$-uniformly smooth Banach spaces. We also define a proximal-point mapping associated with the generalized $\alpha\beta$-$H((.,.),(.,.))$-mixed accretive mapping and show that it is single-valued and Lipschitz continuous. By using the technique of proximal mapping, an iterative
algorithm is constructed in Banach spaces. Under some suitable conditions, we prove the convergence of iterative sequence
generated by the algorithm. The results presented in this paper can
be viewed as an extension and generalization of some known results
\cite{A-D}, \cite{H-G2}-\cite{H-G1}, \cite{K-A-S}-\cite{K-K-S}, \cite{Luo,X-W,Z-H,Z-H1}. For illustration of Definitions 2.8, 3.1,
and Theorem 4.6, Examples 3.2, 3.3 and 4.7 are given, respectively.

\section{Preliminaries}

$\hspace{6mm}$ \noindent Let $X$ be a real Banach space equipped with the norm $\|.\|,$
and let $X^*$ be the topological dual space of $X$. Let
$\langle.,.\rangle$ be the dual pair between $X$ and $X^*$, and
let $2^X$ be the power set of $X$.
\begin{defn}
$\cite{H-X}$ {\rm For $q>1$, a mapping $J_q:X\to 2^{X^*}$ is said to
be a {\it generalized duality mapping}, if it is defined by}
\begin{eqnarray*}
J_q(x)~=~\{f^*\in X^*:~\langle x,f^* \rangle
=\|x\|^q,~\|f^*\|=\|x\|^{q-1}\}, ~~~\forall~x\in X.
\end{eqnarray*}
\end{defn}
\noindent In particular, $J_2$ is the {\it usual normalized duality
mapping} on $X$. It is known that, in general
$$J_q(x)=\|x\|^{q-1}J_2(x)~{\forall}~x(\neq0)~\in X.$$
If $X\equiv H$ a real Hilbert space, then $J_2$ becomes {\it
an identity mapping} on $H.$

\begin{defn}$\cite{H-X}$ {\rm A Banach space $X$ is called {\it smooth} if, for every $x\in X$
with $\|x\|=1,$ there exists a unique $f\in X^*$ such that
$\|f\|=f(x)=1$. \par The {\it modulus of smoothness} of $X$ is a
function $\rho_X:[0,\infty)\to [0,\infty)$, defined by}
\begin{eqnarray*}
\rho_X(t)~=~{\rm
sup}\left\{\frac{1}{2}(\|x+y\|+\|x-y\|)-1:~\|x\|\leq1,~\|y\|\leq
t\right\}.
\end{eqnarray*}
\end{defn}

\begin{defn}$\cite{H-X}$
{\rm A Banach space $X$ is called\\
\noindent(i)~{\it uniformly smooth} if}
\begin{eqnarray*}
&&lim_{t\to 0} \frac{\rho_X(t)}{t}~=~0;
\end{eqnarray*}
\noindent{\rm(ii)~ $q$-{\it uniformly smooth}, for $q>1$, if there
exists a constant $c>0$ such that}
\begin{eqnarray*}
&&\rho_X(t)\leq c~t^q,~t\in[0,\infty).
\end{eqnarray*}
\end{defn}
Note that $J_q$ is single-valued if $X$ is uniformly smooth.
Concerned with the characteristic inequalities in $q$-uniformly
smooth Banach spaces, Xu $\cite{H-X}$ proved the following result.

\begin{lem}
Let $X$ be a real uniformly smooth Banach space. Then $X$ is
$q$-uniformly smooth if and only if there exists a constant $c_q>0$
such that, for all $x,y\in ~X,$
\begin{eqnarray*}
\|x+y\|^q\leq~\|x\|^q+q\langle y,J_q(x)\rangle+c_q \|y\|^q.
\end{eqnarray*}
\end{lem}
From Lemma 2 of Liu \cite{LIU}, it is easy to have the following lemma.

%
In the sequel, we recall important basic concepts and definitions, which will be used in this work.

\begin{defn}
{\rm A mapping $A:X\to X$ is said to be

\noindent(i)~{\it accretive} if
\begin{eqnarray*}
\langle A(x)-A(y),~J_q(x-y)\rangle~\geq~0,~~\forall x,~y\in X;
\end{eqnarray*}

\noindent(ii)~{\it strictly accretive} if
\begin{eqnarray*}
\langle A(x)-A(y),~J_q(x-y)\rangle~>~0,~~\forall x,~y\in X;
\end{eqnarray*}
and equality holds if and only if $x = y$;

\noindent(iii)~ $\xi$-{\it strongly accretive} if there exists a
constant $\xi>0$ such that
\begin{eqnarray*}\langle A(x)-A(y),~J_q(x-y)\rangle~\geq~\xi~\|x-y\|^q,~~\forall
x,~y\in X;
\end{eqnarray*}

\noindent(iv)~ $\mu$-{\it cocoercive} if there exists a constant
$\mu>0$ such that
\begin{eqnarray*}\langle A(x)-A(y),~J_q(x-y)\rangle~\geq~\mu~\|A(x)-A(y)\|^q,~~\forall
x,~y\in X;
\end{eqnarray*}

\noindent(v)~ $\gamma$-{\it relaxed cocoercive} if there exists a
constant $\gamma>0$ such that
\begin{eqnarray*}\langle A(x)-A(y),~J_q(x-y)\rangle~\geq-\gamma~\|A(x)-A(y)\|^q,~~\forall
x,~y\in X;
\end{eqnarray*}

\noindent(vi)~$\zeta$-{\it Lipschitz continuous} if there exists a
constant $\zeta>0$ such that
\begin{eqnarray*}
\|A(x)-A(y)\|~\leq~\zeta~\|x-y\|,~~\forall x,~y\in X;
\end{eqnarray*}

\noindent(vii)~ $\alpha$-{\it expansive} if there exists a constant
$\alpha>0$ such that
\begin{eqnarray*}
\|A(x)-A(y)\|~\geq~\alpha~\|x-y\|,~~\forall x,~y\in X;
\end{eqnarray*}
\par if $\alpha~=~1,$ then it is {\it expansive}.}
\end{defn}

\begin{defn}
$\cite{A-D,Z-H}${\rm Let $H:X\times X\to X$ and $A,B:X\to X$
be the single-valued mappings. Then

\vspace{.25cm}\noindent (i)~$H(A,.)$ is said to be $\alpha$-{\it
strongly accretive} with respect to $A$ if there exists a constant
$\alpha>0$ such that
\begin{eqnarray*}
\langle H(Ax,u)-H(Ay,u),~J_q(x-y)\rangle~\geq~\alpha~\|x-y\|^q,~~\forall
x,~y,~u\in X;
\end{eqnarray*}
(ii)~$H(.,B)$ is said to be $\beta$-{\it relaxed accretive} with
respect to $B$ if there exists a constant $\beta>0$ such that
\begin{eqnarray*}
\langle
H(u,Bx)-H(u,By),~J_q(x-y)\rangle~\geq-\beta~\|x-y\|^q,~~\forall
x,~y,~u\in X;
\end{eqnarray*}
\vspace{.25cm}\noindent (iii)~$H(A,.)$ is said to be $\mu$-{\it
cocoercive} with respect to $A$ if there exists a constant $\mu>0$
such that
\begin{eqnarray*}
\langle H(Ax,~u)-H(Ay,~u),~J_q(x-y)\rangle~\geq~\mu~\|Ax-Ay\|^q,~~\forall
x,~y,~u\in X;
\end{eqnarray*}
(iv)~$H(.,B)$ is said to be $\gamma$-{\it relaxed cocoercive} with
respect to $B$ if there exists a constant $\gamma>0$ such that
\begin{eqnarray*}
\langle
H(u,Bx)-H(u,By),~J_q(x-y)\rangle~\geq-\gamma~\|Bx-By\|^q,~~\forall
x,~y,~u\in X;
\end{eqnarray*}
(v)~$H(A,.)$ is said to be $\tau_1$-{\it Lipschitz continuous} with
respect to $A$ if there exists a constant $\tau_1>0$ such that
\begin{eqnarray*}
\| H(Ax,.)-H(Ay,.)\|~\leq~\tau_1~\|x-y\|,~~\forall x,~y\in X;
\end{eqnarray*}
(vi)~$H(.,B)$ is said to be $\tau_2$-{\it Lipschitz continuous}
with respect to $B$ if there exists a constant $\tau_2>0$ such that
\begin{eqnarray*}
\|H(.,Bx)-H(.,By)\|~\leq~\tau_2~\|x-y\|,~~\forall x,~y\in X.
\end{eqnarray*}}
\end{defn}

\begin{defn}
$\rm\cite{H-G2}${\rm Let $H:(X\times X)\times (X\times X)\to X,$ and $A,B,C,D:X\to
X$ be the single-valued mappings. Then

\vspace{.25cm}\noindent (i)~$H((A,.),(C,.))$ is said to be
$(\mu_1,\gamma_1)$-{\it strongly mixed cocoercive} with respect to $(A,C)$ if
there exist constants $\mu_1,\gamma_1>0$ such that
$$\langle
H((Ax,u),(Cx,u))-H((Ay,u),(Cy,u)),~J_q(x-y)\rangle~\geq~\mu_1~\|Ax-Ay\|^q+\gamma_1~\|x-y\|^q,~\forall
x,~y,~u\in X;$$

\vspace{.25cm}\noindent (ii)~$H((.,B),(.,D))$ is
said to be $(\mu_2,\gamma_2)$-{\it relaxed mixed cocoercive} with
respect to $(B,D)$ if there exist constants $\mu_2,\gamma_2>0$ such
that
$$\langle
H((u,Bx),(u,Dx))-H((u,By),(u,Dy)),~J_q(x-y)\rangle~\geq~-\mu_2~\|Bx-By\|^q+\gamma_2~\|x-y\|^q,~\forall
x,~y,~u\in X;
$$

\vspace{.25cm}\noindent (iii)~$H((A,B),(C,D))$ is
said to be $\mu_1\gamma_1\mu_2\gamma_2$-{\it symmetric mixed cocoercive} with
respect to $(A,C)$ and $(B,D)$ if $H((A,.),(C,.))$ is said to be
$(\mu_1,\gamma_1)$-strongly mixed cocoercive with respect to $(A,C)$ and $H((.,B),(.,D))$ is
said to be $(\mu_2,\gamma_2)$-relaxed mixed cocoercive with
respect to $(B,D)$;

\vspace{.25cm}\noindent (iv)~$H((A,B),(C,D))$ is said to be
$\tau$-{\it mixed Lipschitz continuous} with respect to $A,B,C$ and
$D$ if there exists a constant $\tau>0$ such that
$$
\| H((Ax,Bx),(Cx,Dx))-H((Ay,By),(Cy,Dy))\|~\leq~\tau~\|x-y\|,~~\forall
x,~y\in X.
$$}
\end{defn}


\begin{defn}$\rm{\cite{Luo}}$ {\rm Let $T:X\multimap X$ and $M:X \multimap X$ be the set-valued mapping. Then}

\vspace{.25cm}\noindent{\rm(i)~$T$ is said to be {\it accretive} if
\begin{eqnarray*}
\langle u-v,~J_q(x-y)\rangle~\geq~0~~\forall x,~y\in X,~u\in
Tx,~v\in Ty;
\end{eqnarray*}}
\vspace{.25cm}\noindent{\rm(ii)~$T$ is said to be {\it strictly accretive} if
\begin{eqnarray*}
\langle u-v,~J_q(x-y)\rangle~>~0~~\forall x,~y\in X,~u\in
Tx,~v\in Ty;
\end{eqnarray*}
and equality holds if and only if $x=y$.}

\vspace{.25cm}\noindent {\rm(iii)~$T$ is said to be $\mu'$-{\it strongly accretive}
if there exists a constant $\mu'>0$ such that
\begin{eqnarray*}
\langle u-v,~J_q(x-y)\rangle~\geq~\mu'\|x-y\|^q~~\forall x,~y\in
X,~u\in Tx,~v\in Ty;
\end{eqnarray*}}
\vspace{.25cm}\noindent{\rm(iv)~$T$ is said to be $\gamma'$-{\it relaxed accretive} if
there exists a constant $\gamma'>0$ such that
\begin{eqnarray*}
\langle u-v,~J_q(x-y)\rangle~\geq~-\gamma'\|x-y\|^q~~\forall x,~y\in
X,~u\in Tx,~v\in Ty;
\end{eqnarray*}}
\vspace{.25cm}\noindent{\rm(v)~$M(f,.)$ is said to be $\alpha$-{\it strongly accretive} with respect to $f$
if there exists a constant $\alpha>0$ such that
\begin{eqnarray*}
\langle u-v,~J_q(x-y)\rangle~\geq~\alpha\|x-y\|^q~~\forall x,~y,~w\in
X,~u\in M(f(x),w)~v\in M(f(y),w);
\end{eqnarray*}}
\vspace{.25cm}\noindent{\rm(vi)~$M(.,g)$ is said to be $\beta$-{\it relaxed accretive} with respect to $g$
if there exists a constant $\beta>0$ such that
\begin{eqnarray*}
\langle u-v,~J_q(x-y)\rangle~\geq~-\beta\|x-y\|^q~~\forall x,~y,~w\in
X,~u\in M(w,g(x))~v\in M(w,g(y));
\end{eqnarray*}}
\vspace{.25cm}\noindent{\rm(vii)~$M(.,.)$ is said to be $\alpha\beta$-{\it symmetric accretive} with respect to $f$ and $g$ if $M(f,.)$ is $\alpha$-strongly accretive with respect to $f$ and $M(.,g)$ is $\beta$-relaxed accretive with respect to $g$ with $\alpha \geq \beta$ and $\alpha=\beta$ if and only if $x=y$.}
\end{defn}

\section {Generalized $\alpha\beta$-$H((.,.),(.,.))$-mixed accretive mappings}
~~~~~This section deals with a new concept and properties of
generalized  $\alpha\beta$-$H((.,.),(.,.))$-mixed accretive mappings, which provides a
unifying framework for the existing cocoercive operators, accretive
operators in Banach space and monotone operators in Hilbert space.

\begin{defn} {\rm Let $H:(X\times X)\times (X\times X)\to X$, $f,g:X\to X$
and $A,B,C,D:X\to X$ be single-valued mappings. Let
$H((A,B),(C,D))$ is $\mu_1\gamma_1\mu_2\gamma_2$-symmetric mixed cocoercive with respect
to $(A,C)$ and $(B,D)$. Then the set-valued mapping $M:X\times X\multimap X$ is said
to be generalized $\alpha\beta$-$H((.,.),(.,.))$-{\it mixed accretive} with respect to
$(A,C)$, $(B,D)$ and $(f,g)$ if
\begin{enumerate}
\item[{\rm(i)}] $M$ is $\alpha\beta$-symmetric accretive with respect to $f$ and $g$;
\item[{\rm(ii)}] $(H((.,.),(.,.)) +\rho M(f,g))(X)~=~X$, for all $\rho>0$.
\end{enumerate}}
\end{defn}

\vspace{.25cm}The following example illustrate the Definitions (2.8) and (3.1).
{\begin{example} \rm Let $q=2$ and $X=\mathbb{R}^2$ with usual inner
product defined by $$\langle(x_1,x_2),(y_1,y_2)\rangle=x_1y_1+x_2,y_2.$$ Let $A,B,C,D:\mathbb{R}^2\to \mathbb{R}^2$ be defined
by
\begin{equation*}
Ax
=\binom{4x_1}{4x_2},~Bx=\binom{-3x_1}{-3x_2},~Cx=\binom{2x_1}{2x_2},~Dx=\binom{x_1}{x_2},~\forall~
x=(x_1,x_2)\in \mathbb{R}^2.
\end{equation*}
Suppose that $H:(\mathbb{R}^2\times\mathbb{R}^2)\times(\mathbb{R}^2\times\mathbb{R}^2)\to
\mathbb{R}^2$ is defined by
$$ H((Ax,Bx),(Cx,Dx))~=~Ax+Bx+Cx+Dx.$$ Then $H((A,B),(C,D))$ is symmetric mixed cocoercive with respect to $(A,C)$ and $(B,D)$, and mixed Lipschitz continuous with respect to $A,B,C$ and $D$.
\end{example}}
Indeed, let for any $u\in \mathbb{R}^2$, we have
\begin{eqnarray*}
&&\langle H((Ax,u),(Cx,u))-H((Ay,u),(Cy,u)),~x-y\rangle~\\
&&\hspace{4.5cm}=~\langle Ax+Cx-Ay-Cy,~x-y\rangle\\
&&\hspace{4.5cm}=~\langle Ax-Ay,~x-y\rangle+\langle Cx-Cy,~x-y\rangle\\
&&\hspace{4.5cm}=~\langle(4x_1-4y_1,~4x_2-4y_2),~(x_1-y_1,x_2-y_2)\rangle\\
&&\hspace{5.0cm}~+\langle(2x_1-2y_1,~2x_2-2y_2),~(x_1-y_1,x_2-y_2)\rangle\\
&&\hspace{4.5cm}=~4(x_1-y_1)^2+4(x_2-y_2)^2+2(x_1-y_1)^2~+~2(x_2-y_2)^2\\
&&\hspace{4.5cm}=~4\|x-y\|^2+2\|x-y\|^2
\end{eqnarray*}
and
\begin{eqnarray*}
&&\hspace{2.5cm}\|Ax-Ay\|^2~=~\langle Ax-Ay,~Ax-Ay\rangle\\
&&\hspace{4.5cm}=\langle(4x_1-4x_2,~4y_1-4y_2),~(4x_1-4x_2,~4y_1-4y_2)\rangle\\
&&\hspace{4.5cm}=16(x_1-y_1)^2+16(x_2-y_2)^2,
\end{eqnarray*}
that is,
\begin{eqnarray}
\hspace{-0.9cm}\langle H((Ax,u),(Cx,u))-H((Ay,u),(Cy,u)),~x-y\rangle\geq \frac{1}{4}~\|Ax-Ay\|^2+2\|x-y\|^2.
\end{eqnarray}
Hence, $H((A,B),(C,D))$ is $(\frac{1}{4},2)$-{\it strongly mixed
cocoercive} with respect to $(A,C)$.
\begin{eqnarray*}
&&\langle H((u,Bx),(u,Dx))-H((u,By),(u,Dy)),~x-y\rangle\\
&&\hspace{4.5cm}=~\langle Bx+Dx-By-Dy,~x-y\rangle\\
&&\hspace{4.5cm}=~\langle Bx-By,~x-y\rangle+\langle Dx-Dy,~x-y\rangle\\
&&\hspace{4.5cm}=~\langle(-3x_1+3y_1,~-3x_2+3y_2),~(x_1-y_1,x_2-y_2)\rangle\\
&&\hspace{5.0cm}~+\langle(x_1-y_1,~x_2-y_2),~(x_1-y_1,x_2-y_2)\rangle\\
&&\hspace{4.5cm}=~-3(x_1-y_1)^2-3(x_2-y_2)^2+(x_1-y_1)^2~+~(x_2-y_2)^2\\
&&\hspace{4.5cm}=~-3\|x-y\|^2+\|x-y\|^2
\end{eqnarray*}
and
\begin{eqnarray*}
&&\hspace{2.4cm}\|Bx-By\|^2~=~\langle Bx-By,~Bx-By\rangle\\
&&\hspace{4.5cm}=\langle(-3x_1+3y_1,~-3x_2+3y_2),~(-3x_1+3y_1,~-3x_2+3y_2)\rangle\\
&&\hspace{4.5cm}=-9\{(x_1-y_1)^2+(x_2-y_2)^2\},
\end{eqnarray*}
that is,
\begin{eqnarray}
\langle H((u,Bx),(u,Dx))-H((u,By),(u,Dy)),x-y\rangle\geq -\frac{1}{3}~\|Bx-By\|^2+\|x-y\|^2.
\end{eqnarray}
Hence, $H((A,B),(C,D))$ is $(\frac{1}{3},1)$-{\it relaxed mixed
cocoercive} with respect to $(B,D)$.

\vspace{.25cm}
From (3.1) and (3.2), $H((A,B),(C,D))$ is symmetric mixed cocoercive with respect to $(A,C)$ and $(B,D)$.
\begin{eqnarray*}
\|H((Ax,Bx),(Cx,Dx))-H((Ay,By),(Cy,Dy))\|^2\\
&&\hspace{-5.4cm}~=~\langle
H((Ax,Bx),(Cx,Dx))-H((Ay,By),(Cy,Dy)),\\
&&\hspace{-4.6cm}~H((Ax,Bx),(Cx,Dx))-H((Ay,By),(Cy,Dy))\rangle\\
&&\hspace{-5.4cm}~=~\langle
(Ax+Bx+Cx+Dx)-(Ay+By+Cy+Dy),\\
&&\hspace{-4.6cm}~(Ax+Bx+Cx+Dx)-(Ay+By+Cy+Dy)\rangle\\
&&\hspace{-5.4cm}~=~\langle
(4x_1,4x_2)-(4y_1,4y_2),(4x_1,4x_2)-(4y_1,4y_2)\rangle\\
&&\hspace{-5.4cm}~=~16\{(x_1-y_1)^2+(x_2-y_2)^2\},
\end{eqnarray*}
that is,
$$\hspace{-1.5cm}\|H((Ax,Bx),(Cx,Dx))-H((Ay,By),(Cy,Dy))\|~\leq~4\|x-y\|.$$
Hence, $H((A,B),(C,D))$ is mixed Lipschitz continuous
with respect to $A,B,C$ and $D$.

\vspace{0.3cm} Now, we show that $M(f,g)$ is symmetric accretive with respect to $f$ and $g$.\\
Let $f,g:\mathbb{R}^2\to¨\mathbb{R}^2$ be defined by
\begin{eqnarray*}
f(x)=\binom{5x_1-\frac{2}{3}x_2}{\frac{2}{3}x_1+5x_2},~g(x)=\binom{\frac{7}{4}x_1+\frac{3}{4}x_2}{-\frac{3}{4}x_1+\frac{7}{4}x_2},~\forall~x=(x_1,x_2)\in\mathbb{R}^2.
\end{eqnarray*}
\noindent
Suppose that $M:(\mathbb{R}^2\times\mathbb{R}^2)\to \mathbb{R}^2$ is defined by
\begin{eqnarray*}
M(fx,gx)=fx-gx,~~~~~\forall~~x=(x_1,x_2)\in\mathbb{R}^2.
\end{eqnarray*}
Then $M(f,g)$ is symmetric accretive with respect to $f$ and $g$, and $M$ is generalized $\alpha\beta$-mixed accretive with respect to
$(A,B)$, $(C,D)$ and $(f,g)$.

\vspace{.25cm} Let for any $w\in \mathbb{R}^2$, we have
\begin{eqnarray*}
&&\langle M(fx,w)-M(fy,w),~x-y\rangle~=~\langle fx-w-fy+w,~x-y\rangle\\
&&\hspace{4.5cm}=~\langle fx-fy,~x-y\rangle\\
&&\hspace{4.5cm}=~\langle(5(x_1-y_1)-\frac{2}{3}(x_2-y_2)),(\frac{2}{3}(x_1-y_1)+5(x_2-y_2)),\\
&&\hspace{5.0cm}~(x_1-y_1,x_2-y_2)\rangle\\
&&\hspace{4.5cm}=~5(x_1-y_1)^2+5(x_2-y_2)^2\\
&&\hspace{4.5cm}=~5\|x-y\|^2,
\end{eqnarray*}
that is,
\begin{eqnarray}
&&\hspace{-1.0cm}\langle u-v,~x-y\rangle~\geq~5\|x-y\|^2,~\forall~x,y\in X,~u\in M(fx,w),~v\in M(fy,w).
\end{eqnarray}
Hence, $M(f,g)$ is $5$-strongly accretive with respect to $f$ and
\begin{eqnarray*}
&&\langle M(w,gx)-M(w,gy),~x-y\rangle~=~\langle w-gx-w+gy,~x-y\rangle\\
&&\hspace{4.5cm}=-\langle gx-gy,~x-y\rangle\\
&&\hspace{4.5cm}=-\langle(\frac{7}{4}(x_1-y_1)+\frac{3}{4}(x_2-y_2)),(-\frac{7}{4}(x_1-y_1)+\frac{3}{4}(x_2-y_2)),\\
&&\hspace{5.0cm}~(x_1-y_1,x_2-y_2)\rangle\\
&&\hspace{4.5cm}=-\{\frac{7}{4}(x_1-y_1)^2+\frac{7}{4}(x_2-y_2)^2\}\\
&&\hspace{4.5cm}=-\frac{7}{4}\|x-y\|^2,
\end{eqnarray*}
that is,
\begin{eqnarray}
&&\hspace{-1.0cm}\langle u-v,~x-y\rangle~\geq~-\frac{7}{4}\|x-y\|^2,~\forall~x,y\in X,~u\in M(w,gx),~v\in M(w,gy).
\end{eqnarray}
Hence, $M(f,g)$ is $\frac{7}{4}$-relaxed accretive with respect to $g$.

From (3.3) and (3.4), $M(f,g)$ is symmetric accretive with respect to $f$ and $g$. Also for any $x\in \mathbb{R}^2$, we have
\begin{eqnarray*}
&&\hspace{-1.5cm}[H((A,B),(C,D))+\rho~M(f,g)](x)~=~[H((Ax,Bx),(Cx,Dx))+\rho~M(fx,gx)]\\
&&\hspace{3.3cm}~=~(Ax+Bx+Cx+Dx)+\rho(fx-gx)\\
&&\hspace{3.3cm}~=~(4x_1,4x_2)+(-3x_1,-3x_2)+(2x_1,2x_2)+(x_1,x_2)\\
&&\hspace{3.7cm}~+~\rho\left(\left(5x_1-\frac{2}{3}x_2,\frac{2}{3}x_1+5x_2\right)-\left(\frac{7}{4}x_1-\frac{3}{4}x_2,\frac{3}{4}x_1+\frac{7}{4}x_2\right)\right)\\
&&\hspace{3.3cm}~=~(4x_1,4x_2)+\rho\left(\frac{13}{4}x_1+\frac{1}{12}x_2,-\frac{1}{12}x_1+\frac{13}{4}x_2\right)\\
&&\hspace{3.3cm}~=~\left(\left(4+\frac{13}{4}\rho\right)x_1+\frac{1}{12}\rho x_2,\left(4-\frac{1}{12}\rho\right)x_1+\frac{13}{4}\rho x_2\right),
\end{eqnarray*}
it can be easily verify that the vector on right hand side generate the whole $\mathbb{R}^2.$ Therefore, we have
\begin{eqnarray*}
&&\hspace{-1.5cm}[H((A,B),(C,D))+\rho~M(f,g)]\mathbb{R}^2~=\mathbb{R}^2.
\end{eqnarray*}
Hence, $M$ is generalized $\alpha\beta$-$H((.,.),(.,.))$-mixed accretive with respect to $(A,C)$, $(B,D)$ and $(f,g)$.

{\begin{example} \rm Let $X=l^2$. Then inner product in $l^2$ is defined by
$$\langle x,y\rangle=\sum\limits_{i=1}^{\infty}~=~x_iy_i,~~~~\forall~x,y\in~l^2.$$
Let $A,B,C,D,f,g:l^2\to l^2$ be single-valued mappings are defined by
$$A(x)=-5x-7e_n,~B(x)=5x+5e_n,~C(x)=-3x,~D(x)=2x+3e_n,$$ $${\rm and}~f(x)=2x,~g(x)=x-e_n,~~\forall~
x=(x_1,x_2,.....,x_n,.......)\in~l^2,~{\rm where}~e_n=(0,0,.....,1,.......)\in~l^2.$$\\
Let $H:l^2\times l^2\to l^2$ be a single-valued mapping defined by $$H((Ax,Bx),(Cx,Dx))=Ax+Bx+Cx+Dx,~~~\forall~x\in l^2,$$ and let $M:l^2\times l^2\multimap l^2$ be a set-valued mapping defined by $$M(f(x),g(x))~=~f(x)-g(x).$$ Then
\end{example}
\begin{eqnarray*}
&&\hspace{-1.5cm}\langle M(f(x),w)-M(f(y),w),x-y\rangle~=~\langle f(x)-w-g(y)+w,x-y\rangle\\
&&\hspace{3.2cm}~=~\langle 2x-2y,x-y\rangle\\
&&\hspace{3.2cm}~\geq~2~\|x-y\|^2
\end{eqnarray*}
and
\begin{eqnarray*}
&&\hspace{-1.5cm}\langle M(w,g(x))-M(w,g(y)),x-y\rangle~=~\langle w-g(x)-w+g(y),x-y\rangle\\
&&\hspace{3.2cm}~=~-\langle x-e_n-y+e_n,x-y\rangle\\
&&\hspace{3.2cm}~\geq~-1\|x-y\|^2.
\end{eqnarray*}
Then, for $\rho=1$
\begin{eqnarray*}
&&\hspace{-1.0cm}\|[H((A,B),(C,D))+\rho~M(f,g)](x)\|\\
&&\hspace{3.4cm}=\|Ax+Bx+Cx+Dx+f(x)-g(x)\|\\
&&\hspace{3.3cm}~=~\|-5x-7e_n+5x+5e_n-3x+2x+3e_n+2x-x\|\\
&&\hspace{3.3cm}~\geq~\|e_n\|^2=\langle e_n,e_n\rangle~=~\sum\limits_{i=1}^{\infty}e_{n_i}^2=1
\end{eqnarray*}
and so $0\not\in~[H((A,B),(C,D))+\rho~M(f,g)](l^2)$. Thus, $M$ is not a generalized $\alpha\beta$-$H((A,B),(C,D))$-mixed accretive mapping.

\begin{remark}
\vspace{.20cm} {\rm\noindent(i)~If $H((A,B),(C,D))=H(A,B)$,
then generalized $\alpha\beta$-$H((.,.),(.,.))$-mixed accretive mapping
reduces to generalized $H(.,.)$-accretive mapping considered in
\cite{K-K-S}.}

\vspace{.20cm} {\rm\noindent(ii)~If $H((A,B),(C,D))=B$,
then generalized $\alpha\beta$-$H((.,.),(.,.))$-mixed accretive mapping
reduces to generalized $B$-monotone mapping considered in
\cite{Luo}.}

\vspace{.20cm} {\rm\noindent(iii)~If $H((A,B),(C,D))=H(A,B)$, $M(.,.)=M$ and $M$
is $\eta$-cocoercive, then generalized $\alpha\beta$-$H((.,.),(.,.))$-mixed accretive mapping
reduces to $H(.,.)$-$\eta$-cocoercive mapping considered in
\cite{A-D}.}

%
\vspace{.20cm} {\rm\noindent(iv)~If $H((A,B),(C,D))=H(A,B)$, $M(.,.)=M$ and $M$ is
accretive, then generalized $\alpha\beta$-$H((.,.),(.,.))$-mixed accretive mapping reduces
to $H(.,.)$-accretive mapping considered in \cite{Z-H}.}

\vspace{.20cm} {\rm\noindent(v)~If $H((A,B),(C,D))= H(.)$, $M(.,.)=M$ and $M$ is
accretive (monotone), then generalized $\alpha\beta$-$H((.,.),(.,.))$-mixed accretive mapping reduces to $H$-accretive
mapping considered in \cite{F-H1,F-H2}.}

\vspace{.20cm} {\rm\noindent(vi)~If $X$ is Hilbert space, $M(f,g)=M$ and $M$ is $m$-relaxed monotone,
then generalized $\alpha\beta$-$H((.,.),(.,.))$-mixed accretive mapping
reduces to $H((.,.),(.,.))$-mixed cocoercive mapping considered in
\cite{H-G2}.}
\\
\par {\rm Since generalized $\alpha\beta$-$H((.,.),(.,.))$-mixed accretive mapping is a generalization of the maximal accretive mapping,
it is sensible that there are similar properties between of them.
The following result confirms this expectation.}
\end{remark}

\begin{prop}
Let the set-valued mapping $M :X\multimap X$ be a generalized $\alpha\beta$-$H((.,.),(.,.))$-mixed
accretive mapping with respect to $(A,C)$, $(B,D)$ and $(f,g)$. If $A$ is
$\alpha_1$-expansive, $B$ is $\beta_1$-Lipschitz continuous, and $\alpha>\beta$, $\mu_1>
\mu_2,~\alpha_1>\beta_1$ and $\gamma_1,~\gamma_2>0$. Then
the following inequality:
\begin{eqnarray*}
&&~\langle u-v,~J_q(x-y)\rangle~\geq~0,
\end{eqnarray*}
holds for all~$(v,~y)\in~{\rm Graph}(M(f,g))$,~implies~$(u,~x)\in~M(f,g)$,~where
\begin{eqnarray*}
&&\hspace{1.0cm}{\rm Graph}(M(f,g))~=~\{(u,~x)\in~X\times X:~(u,~x)\in M(f(x),g(x))\}.
\end{eqnarray*}
\end{prop}
\textbf{Proof.} Suppose on the contrary that there exists
$(u_0,~x_0)\not \in {\rm Graph}(M(f,g))$ such that
\begin{eqnarray}
\langle u_0-v,~J_q(x_0-y)\rangle~\geq~0,\forall~(y,v)\in~{\rm Graph}(M(f,g)).
\end{eqnarray}
Since $M$ is generalized $\alpha\beta$-$H((.,.),(.,.))$-mixed accretive with respect to $(A,C)$ and $(B,D)$, we know that
$(H((.,.),(.,.))+\rho~M(f,g))(X)~=~X$ holds for all $\rho>0$ and
so there exists $(u_1,x_1)\in {\rm Graph}(M(f,g))$ such that
\begin{eqnarray}
&&H((Ax_1,Bx_1),(Cx_1,Dx_1)) +\rho u_1=H((Ax_0,Bx_0),(Cx_0,Dx_0))
+\rho u_0\in~X.
\end{eqnarray}
\begin{eqnarray*}
&&\hspace{-1.0cm}{\rm Now,}~~~~\rho~u_0-\rho~u_1~=~H((Ax_1,Bx_1),(Cx_1,Dx_1))-H((Ax_0,Bx_0),(Cx_0,Dx_0))\in~X.\\
&&\langle\rho~u_0-\rho~u_1,~J_q(x_0-x_1)\rangle~=-\langle H((Ax_0,Bx_0),(Cx_0,Dx_0))-H((Ax_1,Bx_1),(Cx_1,Dx_1)),\\
&&\hspace{9.5cm}J_q(x_0-x_1)\rangle.
\end{eqnarray*}
Since $M$ is $\alpha\beta$-symmetric accretive with respect to $f$ and $g$, we obtain
\begin{eqnarray*}
&&\hspace{-0.3cm}(\alpha-\beta)~\|x_0-x_1\|^q\leq\rho~\langle u_0-u_1,~J_q(x_0-x_1)\rangle\\
&&\hspace{2.3cm}=-\langle H((Ax_0,Bx_0),(Cx_0,Dx_0))-H((Ax_1,Bx_1),(Cx_1,Dx_1)),~J_q(x_0-x_1)\rangle\\
&&\hspace{2.3cm}=-\langle H((Ax_0,Bx_0),(Cx_0,Dx_0))-H((Ax_1,Bx_0),(Cx_1,Dx_0)),~J_q(x_0-x_1)\rangle\\
&&\hspace{2.6cm}-\langle H((Ax_1,Bx_0),(Cx_1,Dx_0))-H((Ax_1,Bx_1),(Cx_1,Dx_1)),~J_q(x_0-x_1)\rangle.
\end{eqnarray*}
\vspace{-0.8cm}
\begin{eqnarray}
&&~~~~~~~~
\end{eqnarray}
Since $H((A,B),(C,D))$ is $\mu_1\gamma_1\mu_2\gamma_2$-symmetric mixed cocoercive with respect to $(A,C)$ and $(B,D)$, thus (3.7) becomes
\begin{eqnarray}
&&\hspace{-0.8cm}(\alpha-\beta)~\|x_0-x_1\|^q\leq-\mu_1\| Ax_0-Ax_1\|^q-\gamma_1\|x_0-x_1\|^q+\mu_2\| Bx_0-Bx_1\|^q-\gamma_2\|x_0-x_1\|^q.
\end{eqnarray}
Since $A$ is $\alpha_1$-expansive and $B$ is $\beta_1$-Lipschitz continuous, thus (3.8) becomes
\begin{eqnarray*}
&&\hspace{0.4cm}(\alpha-\beta)~\|x_0-x_1\|^q\leq-\mu_1\alpha_1^q\| x_0-x_1\|^q-\gamma_1\|x_0-x_1\|^q+\mu_2\beta_1^q\| x_0-x_1\|^q-\gamma_2\|x_0-x_1\|^q\\
&&\hspace{3.1cm}=-[(\mu_1\alpha_1^2-\mu_q\beta_1^q)+(\gamma_1+\gamma_2)]~\| x_0-x_1\|^2\\
&&\hspace{2.8cm}0\leq~(\alpha-\beta)~\|x_0-x_1\|^q\leq-[(\mu_1\alpha_1^q-\mu_2\beta_1^q)+(\gamma_1+\gamma_2)]~\| x_0-x_1\|^q\\
&&\hspace{2.8cm}0\leq-(r+m)~\| x_0-x_1\|^q~\leq~0,\\
&&\hspace{-1.1cm}~{\rm where}~r~=~(\mu_1\alpha_1^q-\mu_2\beta_1^q)+(\gamma_1+\gamma_2)~{\rm and}~m~=~(\alpha-\beta),
\end{eqnarray*}
which gives $x_0=x_1$ since $\alpha>\beta,~\mu_1>\mu_2,~\alpha_1>\beta_2,$ and $\gamma_1,\gamma_2>0$. By (3.6), we have $u_0=u_1$, a contradiction. This complete the proof.$\hspace{12.7cm}\Box$

\begin{thm}
Let the set-valued mapping $M :X\multimap X$ be a generalized $\alpha\beta$-$H((.,.),(.,.))$-mixed
accretive mapping with respect to $(A,C)$, $(B,D)$ and $(f,g)$. If $A$ is
$\alpha_1$-expansive, $B$ is $\beta_1$-Lipschitz continuous, and $\alpha>\beta,~\mu_1>
\mu_2,~\alpha_1>\beta_1$ and $\gamma_1,~\gamma_2>0$, then
~$(H((A,B),(C,D))+\rho M(f,g))^{-1}$ is~single-valued.
\end{thm}
\vspace{.3cm} \textbf{Proof.} For any given $x\in X$, let $u, v\in
(H((A,B),(C,D))+\rho M(f,g))^{-1}(x)$. It follows that
\begin{eqnarray*}
\left\{
\begin{array}{ll}
-H((Au,Bu),(Cu,Du))+x\in~\rho M(f,g)u,\\
-H((Av,Bv),(Cv,Dv))+x\in~\rho M(f,g)v.
\end{array}
\right.
\end{eqnarray*}
Since $M$ is $\alpha\beta$-symmetric accretive with respect to $f$ and $g$, we have
\begin{eqnarray*}
&&\hspace{-0.1cm}(\alpha-\beta)\|u-v\|^q~\leq~\frac{1}{\rho}\langle-H((Au,Bu),(Cu,Du))+x-(-H((Av,Bv),(Cv,Dv))+x),~J_q(u-v)\rangle\\
&&\hspace{-0.1cm}(\alpha-\beta)\|u-v\|^q~\leq~\langle-H((Au,Bu),(Cu,Du))+x-(-H((Av,Bv),(Cv,Dv))+x),~J_q(u-v)\rangle\\
&&\hspace{2.3cm}=-\langle H((Au,Bu),(Cu,Du))-H((Av,Bv),(Cv,Dv)),~J_q(u-v)\rangle\\
&&\hspace{2.3cm}=-\langle H((Au,Bu),(Cu,Du))-H((Av,Bu),(Cv,Du)),~J_q(u-v)\rangle\\
&&\hspace{2.5cm}~-\langle
H((Av,Bu),(Cv,Du))-H((Av,Bv),(Cv,Dv)),~J_q(u-v)\rangle.
\end{eqnarray*}
\vspace{-0.8cm}
\begin{eqnarray}
&&~~~
\end{eqnarray}
Since $H((A,B),(C,D))$ is $\mu_1\gamma_1\mu_2\gamma_2$-symmetric mixed cocoercive with respect to $(A,C)$ and $(B,D)$, thus (3.9) becomes
\begin{eqnarray}
&&\hspace{-.50cm}\rho (\alpha-\beta)\|u-v\|^q\leq-\mu_1\| Au-Av\|^q-\gamma_1\|u-v\|^q+\mu_2\|Bu-Bv\|^q-\gamma_2\|u-v\|^q.
\end{eqnarray}
Since $A$ is $\alpha_1$-expansive and $B$ is $\beta_1$-Lipschitz continuous, thus (3.10) becomes
\begin{eqnarray*}
&&\hspace{-.50cm}\rho (\alpha-\beta)\|u-v\|^q\leq-\mu_1\alpha_1^q\| u-v\|^q-\gamma_1\| u-v\|^q+\mu_2\beta_1^q\|u-v\|^q-\gamma_2\| u-v\|^q\\
&&\hspace{2.0cm}=-[(\mu_1\alpha_1^q-\mu_2\beta_1^q)+(\gamma_1+\gamma_2)]~\|u-v\|^q\\
&&\hspace{1.7cm}0\leq~(\alpha-\beta)~\|u-v\|^q\leq-(\mu_1\alpha_1^q-\mu_2\beta_1^q)+(\gamma_1+\gamma_2)~\| u-v\|^q\\
&&\hspace{1.7cm}0\leq-(r+\rho m)~\| u-v\|^q~\leq~0,\\
&&\hspace{-2.1cm}~{\rm where}~r~=~(\mu_1\alpha_1^q-\mu_2\beta_1^q)+(\gamma_1+\gamma_2) ~{\rm and}~m~=~(\alpha-\beta).
\end{eqnarray*}
Since $\alpha>\beta,~\mu_1>\mu_2,~\alpha_1>\beta_2$ and $\gamma_1,\gamma_2>0$, it follows that $\|u-v\|~\leq~0$. This
implies that $u=v$ and so $(H((A,B),(C,D))+\rho M(f,g))^{-1}$ is single-valued.$\hspace{7.5cm}\Box$

\begin{defn}
{\rm Let the set-valued mapping $M :X\multimap X$ be a generalized $\alpha\beta$-
$H((.,.),(.,.))$-mixed accretive mapping with respect to $(A,C)$, $(B,D)$ and $(f,g)$. If $A$ is $\alpha_1$-expansive, $B$ is $\beta_1$-Lipschitz
continuous, and $\alpha>\beta,~\mu_1> \mu_2,~\alpha_1>\beta_1$ and $\gamma_1,~\gamma_2>0$,
then the {\it proximal-point mapping}
$R_{\rho,~M(.,.)}^{H((.,.),(.,.))}:X\to~X$ is defined by
\begin{eqnarray}
&&R_{\rho,~M(.,.)}^{H((.,.),(.,.))}(u)~=~(H((A,B),(C,D))+\rho
M(f,g))^{-1}(u),~~\forall~u\in X.
\end{eqnarray}}
\end{defn}

\begin{remark}

\vspace{.20cm} {\rm\noindent(i)~If $H((A,B),(C,D))=H(A,B)$,
then the proximal-point mapping $R_{\rho,~M(.,.)}^{H((.,.),(.,.))}$ reduces to $R^{H(.,.)}_{M(.,.),\rho}$ considered in
\cite{K-K-S}.}

\vspace{.20cm} {\rm\noindent(ii)~If $H((A,B),(C,D))=B$,
then the proximal-point mapping $R_{\rho,~M(.,.)}^{H((.,.),(.,.))}$ reduces to $R^B_{M(.,.),\rho}$ considered in
\cite{Luo}.}

\vspace{.20cm} {\rm\noindent(iii)~If $H((A,B),(C,D))=H(A,B)$, $M(f,g)=M$, and $M$
is $\eta$-cocoercive, then the  proximal-point mapping $R_{\rho,~M(.,.)}^{H((.,.),(.,.))}$ reduces to $R^{H(.,.)-\eta}_{M,~\rho}$ considered in
\cite{A-D}.}

%
\vspace{.20cm} {\rm\noindent(iv)~If $H((A,B),(C,D))=H(A,B)$, $M(f,g)=M$, and $M$ is
accretive, then the proximal-point mapping $R_{\rho,~M(.,.)}^{H((.,.),(.,.))}$ reduces to $R^{H(.,.)}_{M,~\rho}$ considered in \cite{Z-H}.}

\vspace{.20cm} {\rm\noindent(v)~If $H((A,B),(C,D))= H$, $M(f,g)=M$, and $M$ is
accretive (monotone), then the proximal-point mapping $R_{\rho,~M(.,.)}^{H((.,.),(.,.))}$ reduces to $R^{H}_{M,~\rho}$ considered in \cite{F-H1,F-H2}.}

\vspace{.20cm} {\rm\noindent(vi)~If $X$ is Hilbert space, $M(f,g)=M$, and $M$ is $m$-relaxed monotone,
then the proximal-point mapping $R_{\rho,~M(.,.)}^{H((.,.),(.,.))}$ reduces to the resolvent operator $R^{H((.,.),(.,.))}_{\rho,~M(.,.)}$ considered in
\cite{H-G2}.}
\end{remark}
\par Now we prove that the proximal-point mapping
defined by (3.11) is Lipschitz continuous.
\begin{thm}
Let the set-valued mapping $M :X\multimap X$ be a generalized $\alpha\beta$-
$H((.,.),(.,.))$-mixed accretive mapping with respect to $(A,C)$, $(B,D)$ and $(f,g)$. If $A$ is
$\alpha_1$-expansive, $B$ is $\beta_1$-Lipschitz continuous, and $\alpha>\beta,~\mu_1>
\mu_2,~\alpha_1>\beta_1$ and $\gamma_1,~\gamma_2>0$, then
the proximal-point mapping $R_{\rho,~M(.,.)}^{H((.,.),(.,.))}:X\to~X$ is
$\frac{1}{r+\rho m}$-Lipschitz continuous, that is,
\begin{eqnarray*}
\|R_{\rho,~M(.,.)}^{H((.,.),(.,.))}(u)-R_{\rho,~M(.,.)}^{H((.,.),(.,.))}(v)\|\leq\frac{1}{r+\rho m}\|u-v\|,~\forall~u,v\in
X.
\end{eqnarray*}
\end{thm}
\textbf{Proof.} Let $u,~v\in X$ be any given points, It follows from
(3.11) that
\begin{eqnarray*}
R_{\rho,~M(.,.)}^{H((.,.),(.,.))}(u)~=~(H((A,B),(C,D))+\rho M(f,g))^{-1}(u),\\
R_{\rho,~M(.,.)}^{H((.,.),(.,.))}(v)~=~(H((A,B),(C,D))+\rho M(f,g))^{-1}(v).
\end{eqnarray*}
\begin{eqnarray*}
\left\{
\begin{array}{ll}
\frac{1}{\rho}\Big(u-H\Big(\Big(A\Big(R_{\rho,~M(.,.)}^{H((.,.)(.,.))}(u)\Big),B\Big(R_{\rho,~M}^{H((.,.)(.,.))}(u)\Big)\Big),\\
\hspace{1.7cm}\Big(C\Big(R_{\rho,~M(.,.)}^{H((.,.),(.,.))}(u)\Big),~D\Big(R_{\rho,~M}^{H((.,.)(.,.))}(u)\Big)\Big)\Big)\Big)\in~M~\Big( f\Big(R_{\rho,~M(.,.)}^{H((.,.),(.,.))}(u)\Big),g\Big(R_{\rho,~M(.,.)}^{H((.,.),(.,.))}(u)\Big)\Big)\\
\frac{1}{\rho}\Big(v-H\Big(\Big(A\Big(R_{\rho,~M(.,.)}^{H((.,.)(.,.))}(v)\Big),B\Big(R_{\rho,~M}^{H((.,.)(.,.))}(v)\Big)\Big),\\
\hspace{1.7cm}\Big(C\Big(R_{\rho,~M(.,.)}^{H((.,.),(.,.))}(v)\Big),~D\Big(R_{\rho,~M}^{H((.,.)(.,.))}(v)\Big)\Big)\Big)\Big)\in~M~\Big( f\Big(R_{\rho,~M(.,.)}^{H((.,.),(.,.))}(v)\Big),g\Big(R_{\rho,~M(.,.)}^{H((.,.),(.,.))}(v)\Big)\Big).\\
\end{array}
\right.
\end{eqnarray*}
Since $M$ is $\alpha\beta$-symmetric accretive with respect to $f$ and $g$, we have
\begin{eqnarray*}
\left\{
\begin{array}{ll}
\langle\frac{1}{\rho}(u-H((A(R_{\rho,~M(.,.)}^{H((.,.)(.,.))}(u)),B(R_{\rho,~M(.,.)}^{H((.,.)(.,.))}(u))),(C(R_{\rho,~M(.,.)}^{H((.,.)(.,.))}(u)),~D(R_{\rho,~M(.,.)}^{H((.,.)(.,.))}(u))))\\
\hspace{2.2cm}-(v-H((A(R_{\rho,~M(.,.)}^{H((.,.)(.,.))}(v)),B(R_{\rho,~M(.,.)}^{H((.,.)(.,.))}(v))),(C(R_{\rho,~M(.,.)}^{H((.,.)(.,.))}(v)),~D(R_{\rho,~M(.,.)}^{H((.,.)(.,.))}(v)))))),\\
\hspace{3.2cm}~J_q(R_{\rho,~M(.,.)}^{H((.,.)(.,.))}(u)-R_{\rho,~M(.,.)}^{H((.,.)(.,.))}(v))\rangle\geq~(\alpha-\beta)~\|R_{\rho,~M(.,.)}^{H((.,.)(.,.))}(u)-R_{\rho,~M(.,.)}^{H((.,.)(.,.))}(v)\|^q,\\~\\
\langle\frac{1}{\rho}(u-v-H((A(R_{\rho,~M(.,.)}^{H((.,.)(.,.))}(u)),B(R_{\rho,~M(.,.)}^{H((.,.)(.,.))}(u))),(C(R_{\rho,~M(.,.)}^{H((.,.)(.,.))}(u)),~D(R_{\rho,~M(.,.)}^{H((.,.)(.,.))}(u)))))\\
\hspace{2.2cm}+H((A(R_{\rho,~M(.,.)}^{H((.,.)(.,.))}(v)),B(R_{\rho,~M(.,.)}^{H((.,.)(.,.))}(v))),(C(R_{\rho,~M(.,.)}^{H((.,.)(.,.))}(v)),~D(R_{\rho,~M(.,.)}^{H((.,.)(.,.))}(v)))),\\
\hspace{3.2cm}~J_q(R_{\rho,~M(.,.)}^{H((.,.)(.,.))}(u)-R_{\rho,~M(.,.)}^{H((.,.)(.,.))}(v))\rangle\geq(\alpha-\beta)~\|R_{\rho,~M(.,.)}^{H((.,.)(.,.))}(u)-R_{\rho,~M(.,.)}^{H((.,.)(.,.))}(v)\|^q,
\end{array}
\right.
\end{eqnarray*}
which implies
\begin{eqnarray*}
&&\hspace{-0.3cm}\langle u-v,~J_q(R_{\rho,~M(.,.)}^{H((.,.)(.,.))}(u)-R_{\rho,~M(.,.)}^{H((.,.)(.,.))}(v))\rangle~\geq\\
&&\hspace{1.3cm}~\langle H((A(R_{\rho,~M(.,.)}^{H((.,.)(.,.))}(u)),B(R_{\rho,~M(.,.)}^{H((.,.)(.,.))}(u))),(C(R_{\rho,~M(.,.)}^{H((.,.)(.,.))}(u)),D(R_{\rho,~M(.,.)}^{H((.,.)(.,.))}(u))))\\
&&\hspace{1.5cm}-H((A(R_{\rho,~M(.,.)}^{H((.,.)(.,.))}(v)),B(R_{\rho,~M(.,.)}^{H((.,.)(.,.))}(v))),(C(R_{\rho,~M(.,.)}^{H((.,.)(.,.))}(v)),~D(R_{\rho,~M(.,.)}^{H((.,.)(.,.))}(v)))),\\
&&\hspace{2.2cm}~J_q(R_{\rho,~M(.,.)}^{H((.,.)(.,.))}(u)-R_{\rho,~M(.,.)}^{H((.,.)(.,.))}(v))\rangle+\rho
(\alpha-\beta)~\|R_{\rho,~M(.,.)}^{H((.,.)(.,.))}(u)-R_{\rho,~M(.,.)}^{H((.,.)(.,.))}(v)\|^q.
\end{eqnarray*}
Now, we have
\begin{eqnarray*}
&&\hspace{-0.8cm}~\|u-v\|~\|R_{\rho,~M(.,.)}^{H((.,.)(.,.))}(u)-R_{\rho,~M(.,.)}^{H((.,.)(.,.))}(v)\|^{q-1}~\geq~\langle u-v,~R_{\rho,~M(.,.)}^{H((.,.)(.,.))}(u)-R_{\rho,~M(.,.)}^{H((.,.)(.,.))}(v)\rangle\\
&&\hspace{3.1cm}\geq~\langle H((A(R_{\rho,~M(.,.)}^{H((.,.)(.,.))}(u)),B(R_{\rho,~M(.,.)}^{H((.,.)(.,.))}(u))),(C(R_{\rho,~M(.,.)}^{H((.,.)(.,.))}(u)),D(R_{\rho,~M(.,.)}^{H((.,.)(.,.))}(u))))\\
&&\hspace{3.6cm}-H((A(R_{\rho,~M(.,.)}^{H((.,.)(.,.))}(v)),B(R_{\rho,~M(.,.)}^{H((.,.)(.,.))}(v))),(C(R_{\rho,~M(.,.)}^{H((.,.)(.,.))}(v)),D(R_{\rho,~M(.,.)}^{H((.,.)(.,.))}(v)))),\\
&&\hspace{3.9cm}~J_q(R_{\rho,~M(.,.)}^{H((.,.)(.,.))}(u)-R_{\rho,~M(.,.)}^{H((.,.)(.,.))}(v))\rangle+\rho
(\alpha-\beta)~\|R_{\rho,~M(.,.)}^{H((.,.)(.,.))}(u)-R_{\rho,~M(.,.)}^{H((.,.)(.,.))}(v)\|^q\\
&&\hspace{3.1cm}=~\langle H((A(R_{\rho,~M(.,.)}^{H((.,.)(.,.))}(u)),B(R_{\rho,~M(.,.)}^{H((.,.)(.,.))}(u))),(C(R_{\rho,~M(.,.)}^{H((.,.)(.,.))}(u)),D(R_{\rho,~M(.,.)}^{H((.,.)(.,.))}(u))))\\
&&\hspace{3.6cm}-H((A(R_{\rho,~M(.,.)}^{H((.,.)(.,.))}(v)),B(R_{\rho,~M(.,.)}^{H((.,.)(.,.))}(u))),(C(R_{\rho,~M(.,.)}^{H((.,.)(.,.))}(v)),D(R_{\rho,~M(.,.)}^{H((.,.)(.,.))}(u)))),\\
&&\hspace{5.9cm}~J_q(R_{\rho,~M(.,.)}^{H((.,.)(.,.))}(u)-R_{\rho,~M(.,.)}^{H((.,.)(.,.))}(v))\rangle\\
&&\hspace{3.6cm}+~\langle H((A(R_{\rho,~M(.,.)}^{H((.,.)(.,.))}(v)),B(R_{\rho,~M(.,.)}^{H((.,.)(.,.))}(u))),(C(R_{\rho,~M(.,.)}^{H((.,.)(.,.))}(v)),D(R_{\rho,~M(.,.)}^{H((.,.)(.,.))}(u))))\\
&&\hspace{3.6cm}-H((A(R_{\rho,~M(.,.)}^{H((.,.)(.,.))}(v)),B(R_{\rho,~M(.,.)}^{H((.,.)(.,.))}(v))),(C(R_{\rho,~M(.,.)}^{H((.,.)(.,.))}(v)),D(R_{\rho,~M(.,.)}^{H((.,.)(.,.))}(v)))),\\
&&\hspace{3.9cm}~J_q(R_{\rho,~M(.,.)}^{H((.,.)(.,.))}(u)-R_{\rho,~M(.,.)}^{H((.,.)(.,.))}(v))\rangle+\rho
(\alpha-\beta)~\|R_{\rho,~M(.,.)}^{H((.,.)(.,.))}(u)-R_{\rho,~M(.,.)}^{H((.,.)(.,.))}(v)\|^q.
\end{eqnarray*}
Since $H((A,B),(C,D))$ is $\mu_1\gamma_1\mu_2\gamma_1$-symmetric mixed cocoercive with respect to $(A,C)$ and $(B,D)$, we have
\begin{eqnarray*}
&&\hspace{-.6cm}~\|u-v\|~\|R_{\rho,~M(.,.)}^{H((.,.)(.,.))}(u)-R_{\rho,~M(.,.)}^{H((.,.)(.,.))}(v)\|^{q-1}\\
&&\hspace{3.1cm}\geq~\mu_1\|A(R_{\rho,~M(.,.)}^{H((.,.)(.,.))}(u))-A(R_{\rho,~M(.,.)}^{H((.,.)(.,.))}(v))\|^q+\gamma_1\|R_{\rho,~M(.,.)}^{H((.,.)(.,.))}(u)-R_{\rho,~M(.,.)}^{H((.,.)(.,.))}(v)\|^q\\
&&\hspace{3.6cm}-\mu_2\|B(R_{\rho,~M(.,.)}^{H((.,.)(.,.))}(u))-B(R_{\rho,~M(.,.)}^{H((.,.)(.,.))}(v))\|^q+\gamma_2\|R_{\rho,~M(.,.)}^{H((.,.)(.,.))}(u)-R_{\rho,~M(.,.)}^{H((.,.)(.,.))}(v)\|^q\\
&&\hspace{3.6cm}+\rho (\alpha-\beta)~\|R_{\rho,~M(.,.)}^{H((.,.)(.,.))}(u)-R_{\rho,~M(.,.)}^{H((.,.)(.,.))}(v)\|^q.
\end{eqnarray*}
Since $A$ is $\alpha_1$-expansive and $B$ is $\beta_1$-Lipschitz continuous, we have
\begin{eqnarray*}
&&\hspace{-0.7cm}\|u-v\|~\|R_{\rho,~M(.,.)}^{H((.,.)(.,.))}(u)-R_{\rho,~M(.,.)}^{H((.,.)(.,.))}(v)\|^{q-1}\geq[(\mu_1\alpha_1^q-\mu_2\beta_1^q)+(\gamma_1+\gamma_2)]~\|R_{\rho,~M(.,.)}^{H((.,.)(.,.))}(u)-R_{\rho,~M(.,.)}^{H((.,.)(.,.))}(v)\|^q\\
&&\hspace{5.5cm}+\rho (\alpha-\beta)~\|R_{\rho,~M(.,.)}^{H((.,.)(.,.))}(u)-R_{\rho,~M(.,.)}^{H((.,.)(.,.))}(v)\|^q\\
&&\hspace{5.2cm}\geq~(r+\rho m)~\|R_{\rho,~M(.,.)}^{H((.,.)(.,.))}(u)-R_{\rho,~M(.,.)}^{H((.,.)(.,.))}(v)\|^q,
\end{eqnarray*}
where~$r~=~(\mu_1\alpha_1^q-\mu_2\beta_1^q)+(\gamma_1+\gamma_2)~{\rm and}~m=(\alpha-\beta)$.\\
Hence,
\begin{eqnarray*}
&&\hspace{-1.0cm}\|u-v\|~\|R_{\rho,~M(.,.)}^{H((.,.)(.,.))}(u)-R_{\rho,~M(.,.)}^{H((.,.)(.,.))}(v)\|^{q-1}~\geq~(r+\rho m)~\|R_{\rho,~M(.,.)}^{H((.,.)(.,.))}(u)-R_{\rho,~M(.,.)}^{H((.,.)(.,.))}(v)\|^q,
\end{eqnarray*}
that is
\begin{eqnarray*}
&&\|R_{\rho,~M(.,.)}^{H((.,.),(.,.))}(u)-R_{\rho,~M(.,.)}^{H((.,.),(.,.))}(v)\|\leq\frac{1}{r+\rho m}~\|u-v\|,~\forall~u,v\in
X.
\end{eqnarray*}
This completes the proof.$\hspace{10.0cm}\Box$

\section {An application of generalized $\alpha\beta$-$H((.,.),(.,.))$-mixed accretive mappings for solving variational inclusions.}

\par~~~~~ In this section, we shall show that under suitable assumptions, the generalized $\alpha\beta$-$H((.,.),(.,.))$-mixed accretive mapping
can also play important roles for solving the generalized set-valued variational inclusion in Banach space.

\vspace{.25cm}\noindent Let $S,T:X\multimap CB(X)$ be the set-valued
mappings, and let $f,g:X\to X,$ $A,B,C,D:X\to X$, $F:X\times X\to X$
and $H:(X\times X)\times (X\times X)\to X$ be single-valued mappings.
Suppose that $M:X\times X\multimap X$ is a set-valued mapping such that $M$ be a generalized $\alpha\beta$-
$H((.,.),(.,.))$-mixed accretive mapping with respect to $(A,C)$, $(B,D)$
and $(f,g)$. We consider the following generalized set-valued variational inclusion: for given
$\omega\in X,$ find $u\in X$, $v\in S(u)$ and $w\in T(u)$ such that
\begin{eqnarray}
\omega\in~F(v,w)+M(f(u),g(u)).
\end{eqnarray}

\vspace{0.25cm}\noindent
If $S,T:X\to X$ be single-valued mappings and $M(.,.)=\lambda N(.)$, where $\rho>0$ is a constant, then the problem (4.1)
reduces to the following problem: find $u\in X$ such that
\begin{eqnarray}
\omega\in~F(S(u),T(u))+\lambda N(u).
\end{eqnarray}

\vspace{0.25cm}\noindent
If $M$ is an $(A,\eta)$-accretive mapping, then the problem (4.2) was introduced and studied by Lan et al. \cite{L-C-V}.

\vspace{0.25cm}\noindent If $\lambda=1$, $a=0$ and $F(S(u),T(u))=T(u)$ for all $u\in X$, where $T:X\to X$ is a single-valued mapping, then the problem
(4.2) reduces to the following problem: find $u\in X$ such that
\begin{eqnarray}
0\in~T(u)+N(u).
\end{eqnarray}

\vspace{0.25cm}\noindent
If $N$ is an $H(.,.)$-accretive mapping, then the problem (4.3) was studied by Zou and Huang \cite{Z-H}; and $N$ is a generalized
$m$-accretive mapping, then the problem (4.3) was studied by Bi et al. \cite{B-H-F}.

\vspace{0.25cm}\noindent If $X=H$ is a Hilbert space and $N$ is an $H$-monotone mappings, then the problem (4.3) was introduced and
studied by Fang and Huang \cite{F-H1} and includes many variational inequalities (inclusions) and complementarity problems as
special cases. For example, see \cite{Peng1,Peng2}.

\begin{lem} Let $S,T:X\multimap CB(X)$ be the set-valued
mappings, and let $f,g:X\to X,$ $A,B,C,D:X\to X$, $F:X\times X\to X$
and $H:(X\times X)\times (X\times X)\to X$ be single-valued mappings.
Suppose that $M:X\times X\multimap X$ is a set-valued mapping such that $M$ be a generalized $\alpha\beta$-
$H((.,.),(.,.))$-mixed accretive mapping with respect to $(A,C)$, $(B,D)$
and $(f,g)$. Then $u\in X$, $v\in S(u)$ and $w\in T(u)$ is a solution of problem (4.1) if and only if $u\in X$, $v\in S(u)$ and $w\in T(u)$ satisfies the following relation:
\begin{eqnarray}
u~=~R^{H((.,.),(.,.))}_{\rho,M(.,.)}~[H((Au,Bu),(Cu,Du))-\rho F(v,w)+\rho\omega],
\end{eqnarray}
where $\rho~>~0$ is a constant and $R^{H((.,.),(.,.))}_{\rho,M(.,.)}$ is
the proximal-point mapping defined by (3.11).
\end{lem}

\textbf{Proof.} Observe that for $\rho>0,$
\begin{eqnarray*}
&&\omega\in~F(w,v)+ M(f(u),g(u))\\
&&\hspace{-0.8cm} \Leftrightarrow [H((Au,Bu),(Cu,Du))-\rho F(v,w)+\rho\omega]\in H((Au,Bu),(Cu,Du))+\rho M(f(u),g(u))\\
&&\hspace{-0.8cm} \Leftrightarrow [H((Au,Bu),(Cu,Du))-\rho F(v,w)+\rho\omega]\in (H((A,B),(C,D))+\rho M(f,g))u\\
&&\hspace{-0.8cm} \Leftrightarrow u=(H((A,B),(C,D))+\rho M(f,g))^{-1}[H((Au,Bu),(Cu,Du))-\rho F(v,w)+\rho\omega]\\
&&\hspace{-0.8cm} \Leftrightarrow u=R^{(H(.,.),(.,.))}_{\rho~M(.,.)}[H((Au,Bu),(Cu,Du))-\rho F(v,w)+\rho\omega].\hspace{5cm}\Box
\end{eqnarray*}

\vspace{0.2cm}
\begin{remark}
{\rm We can rewrite the equality (4.4) as:
\begin{eqnarray*}
&&z~=~H((Au,Bu),(Cu,Du))-\rho F(v,w)+\rho\omega,~~~ u=R^{H((.,.),(.,.))}_{\rho,~M(.,.)}(z).
\end{eqnarray*}
where $\omega\in X$ is any given element and $\rho> 0$ is a constant. By Nadler \cite{Nad}, we know that this fixed point formulation enables
us to suggest the following iterative algorithm.}
\end{remark}
\noindent
\begin{alg}
{\rm For any given $z_0\in X$, we can choose $u_0\in X$ such that sequences $\{u_n\}$, $\{v_n\}$ and $\{w_n\}$ satisfy}
\begin{eqnarray*}
\left\{
\begin{array}{ll}
u_n=R^{H((.,.),(.,.))}_{\rho,~M(.,.)}(z_n),\\
v_n\in S(u_n),~~\|v_n-v_{n+1}~\|\leq~\left(1+\frac{1}{n+1}\right)~{\cal{D}}(S(u_n),S(u_{n+1})),\\
w_n\in T(u_n),~\|w_n-w_{n+1}\|~\leq~\left(1+\frac{1}{n+1}\right)~{\cal{D}}(T(u_n),T(u_{n+1})),\\
z_{n+1}~=~H((Au_n,Bu_n),(Cu_n,Du_n))-\rho F(v_n,w_n)+\rho\omega+e_n,\\
\sum\limits_{j=1}^{\infty}\|e_{j}-e_{j-1}\|~\varpi^{-j}~<~\infty,~\forall~\varpi\in~(0,1),~{\rm lim}_{n\to \infty}e_n=0,
\end{array}
\right.
\end{eqnarray*}
{\rm where $\rho > 0$ is a constant, $\omega\in X$ is any given element and $e_n\subset X$ is an error to take into account a possible inexact
computation of the proximal-point mapping point for all $n\geq0$, and ${\cal{D}}(.,.)$ is the Hausdorff metric on} ${\rm CB}(X)$.
\end{alg}

We need the following definitions which will be used to state and prove the main result.
\begin{defn} {\rm A set-valued mapping $G: X\multimap CB(X)$ is said to be $\cal{D}$-{\it Lipschitz continuous} if there exists a constant $l>0$ such
that}
\begin{eqnarray*}
{\cal{D}}(Gx,Gy)\leq~l~\|x-y\|,~~\forall x,~y\in X.
\end{eqnarray*}
\end{defn}

\begin{defn} {\rm Let $S,T:X\multimap X$ be the set-valued mappings, $A,B,C,D:X\to X,~F: X\times X\to X$ and $H: (X\times X)\times(X\times X)\to X$ be single-valued mappings. Then}

\vspace{.25cm}
\noindent{\rm(i)~$F$ is said to be $\sigma$-{\it strongly accretive} with respect to $S$ and $H((A,B),(C,D))$ in the first argument if there exists a constant $\sigma>0$ such that}
\begin{eqnarray*}
&&\hspace{-1.0cm}\langle F(v_1,.)-F(v_{2},.),~J_q(H((Au,Bu),(Cu,Du))-H((Av,Bv),(Cv,Dv)))\rangle~\\
&&\hspace{4.0cm}\geq~\sigma~\|H((Au,Bu),(Cu,Du))-H((Av,Bv),(Cv,Dv))\|^q,\\
&&\hspace{4.5cm}~~\forall~u,v\in~X ~{\rm and}~v_1\in S(u),~v_2\in S(v);
\end{eqnarray*}

\vspace{.25cm}
\noindent{\rm(ii)~$F$ is said to be $\delta$-{\it strongly accretive} with respect to $T$ and $H((A,B),(C,D))$ in the second argument if there exists a constant $\delta>0$ such that}
\begin{eqnarray*}
&&\hspace{-1.0cm}\langle F(.,w_1)-F(.,w_{2}),~J_q(H((Au,Bu),(Cu,Du))-H((Av,Bv),(Cv,Dv)))\rangle~\\
&&\hspace{4.0cm}\geq~\delta~\|H((Au,Bu),(Cu,Du))-H((Av,Bv),(Cv,Dv))\|^q,\\
&&\hspace{4.5cm}~~\forall~u,v\in~X ~{\rm and}~w_1\in T(u),~w_2\in T(v);
\end{eqnarray*}

\vspace{.25cm}
\noindent{\rm(iii)~$F$ is said to be $\epsilon_1$-{\it
Lipschitz continuous} in the first argument if there exists a
constant $\epsilon_1>0$ such that}
\begin{eqnarray*}
\|F(u,v')-F(v,v')\|~\leq~\epsilon_1~\|u-v\|,~~\forall u,~v,~v'\in X;
\end{eqnarray*}

\vspace{.25cm}
\noindent{\rm(iv)~$F$ is said to be $\epsilon_2$-{\it
Lipschitz continuous} in the second argument if there exists a
constant $\epsilon_2>0$ such that}
\begin{eqnarray*}
\|F(v',u)-F(v',v)\|~\leq~\delta_2~\|u-v\|,~~\forall u,~v,~v'\in X.
\end{eqnarray*}
\end{defn}

\par Next, we find the convergence of iterative algorithm for generalized set-valued variational inclusion (4.1).

\begin{thm}
Let $S,T:X\multimap CB(X)$ be the set-valued
mappings, and let $f,g:X\to X,$ $A,B,C,D:X\to X$, $F:X\times X\to X$
and $H:(X\times X)\times (X\times X)\to X$ be single-valued mappings.
Suppose that $M:X\times X\multimap X$ is a set-valued mapping such that $M$ be a generalized $\alpha\beta$-
$H((.,.),(.,.))$-mixed accretive mapping with respect to $(A,C)$, $(B,D)$
and $(f,g)$. Assume that

\vspace{.2cm} \noindent(i) $S$ and $T$ are ${\cal D}$-Lipschitz
continuous with constants $l_1$ and $l_2$, respectively;

\vspace{.2cm} \noindent(ii) $A$ is $\alpha_1$-expansive and $B$ is
$\beta_1$-Lipschitz continuous;

\vspace{.2cm} \noindent(iii) $H((A,B),(C,D))$ is $\tau$-mixed Lipschitz continuous with respect to $A,B,C$ and $D$;

\vspace{.2cm} \noindent(iv) $F$ is is $\sigma$-strongly accretive with respect to $S$ and $H((A,B),(C,D))$ in the first argument and $\delta$-strongly accretive with respect to $T$ and $H((A,B),(C,D))$ in the second argument;

\vspace{.2cm} \noindent(v) $F$ is is $\epsilon_1$-Lipschitz
continuous in the first argument and $\epsilon_2$-Lipschitz
continuous in the second argument;

\vspace{0.2cm}\noindent(vi)~$0<\sqrt[q]{\tau^q+c_q \rho^q(\epsilon_1 l_1+\epsilon_2 l_2)^q-\rho q (\sigma+\delta)\tau^q}<{r+\rho m};$
\vspace{-0.8cm}
\begin{eqnarray}
&&~~
\end{eqnarray}
~where~$r~=~(\mu_1\alpha_1^q~-~\mu_2\beta_1^q)+(\gamma_1+\gamma_2)$ and $m~=~\alpha-\beta$,
and $\alpha>\beta$, $\mu_1>\mu_2,~\alpha_1>\beta_1$ and $\gamma_1,~\gamma_2,\rho>0~$.

\vspace{0.2cm} Then generalized set-valued variational inclusion problem
(4.1) has a solution $(u,v,w)$, where $u\in X,$ $v\in S(u)$ and
$w\in T(u)$, and the iterative sequences
$\{u_n\},~\{v_n\}$ and $\{w_n\}$, generated by Algorithms 4.3
converges strongly to $u,~v$ and $w$, respectively.
\end{thm}

\textbf{Proof.} Since $S$ and $T$ are ${\cal D}$-Lipschitz continuous
with constants $l_1$ and $l_2$, respectively, it follows from
Algorithms 4.3 such that
\begin{eqnarray}
&&\hspace{-1cm}\|v_{n+1}-v_{n}\|\leq
\left(1+\frac{1}{n+1}\right)~{\cal{D}}(S(u_{n+1}),S(u_{n}))\leq
\left(1+\frac{1}{n+1}\right)~l_1~\|u_{n+1}-u_{n}\|,\\
&&\hspace{-1cm}\|w_{n+1}-w_{n}\|\leq
\left(1+\frac{1}{n+1}\right)~{\cal{D}}(T(u_{n+1}),T(u_{n}))\leq
\left(1+\frac{1}{n+1}\right)~l_2~\|u_{n+1}-u_{n}\|,
\end{eqnarray}
for $n=0,~1,~2,~....$

\noindent It follows from (4.4) and Theorem 3.8 that
\begin{eqnarray}
&&\hspace{-1cm}\|u_{n+1}-u_{n}\|~\leq~\|R^{H((.,.),(.,.))}_{\rho,~M(.,.)}(z_{n+1})-R^{H((.,.),(.,.))}_{\rho,~M(.,.)}(z_{n})\|~=~\frac{1}{r+\rho m}\|z_{n+1}-z_{n}\|.
\end{eqnarray}

\noindent Now, we estimate $\|z_{n+1}-z_n\|$ by using
Algorithms 4.3, we have
\begin{eqnarray*}
&&\hspace{-2.4cm}\|z_{n+1}-z_n\|~=~\|[H((Au_{n},Bu_{n}),(Cu_{n},u_{n}))-\rho F(v_n,w_n)+\rho\omega+e_n]\\
&&\hspace{-0.2cm}~-[H((Au_{n-1},Bu_{n-1}),(Cu_{n-1},u_{n-1}))-\rho F(v_{n-1},w_{n-1})+\rho\omega+e_{n-1}]\|\\
&&\hspace{-0.8cm}~\leq~\|H((Au_{n},Bu_{n}),(Cu_{n},u_{n}))-H((Au_{n-1},Bu_{n-1}),(Cu_{n-1},u_{n-1}))
\end{eqnarray*}
\vspace{-0.7cm}
\begin{eqnarray}
&&\hspace{-3.1cm}~+~(\rho F(v_n,w_n)-\rho F(v_{n-1},w_{n-1})\|~+~\|e_n-e_{n-1}\|.
\end{eqnarray}
\noindent By Lemma 2.4, we have
\begin{eqnarray*}
&&\hspace{-0.6cm}\|H((Au_{n},Bu_{n}),(Cu_{n},Du_{n}))-H((Au_{n-1},Bu_{n-1}),(Cu_{n-1},Du_{n-1}))-\rho(F(v_n,w_n)-F(v_{n-1},w_{n-1})\|^q\\
&&\hspace{-0.4cm}~\leq~\|H((Au_{n},Bu_{n}),(Cu_{n},Du_{n}))-H((Au_{n-1},Bu_{n-1}),(Cu_{n-1},Du_{n-1}))\|^q+c_q\rho^q\|F(v_n,w_n)-F(v_{n-1},w_{n-1})\|^q\\
&&\hspace{0.2cm}-~\rho q\langle F(v_n,w_n)-F(v_{n-1},w_{n-1}),~J_q(H((Au_{n},Bu_{n}),(Cu_{n},Du_{n}))-H((Au_{n-1},Bu_{n-1}),(Cu_{n-1},Du_{n-1})))\rangle.
\end{eqnarray*}
\vspace{-0.7cm}
\begin{eqnarray}
~~
\end{eqnarray}
From (iii), we get
\begin{eqnarray}
&&\hspace{-0.7cm}\|H((Au_{n},Bu_{n}),(Cu_{n},Du_{n}))-H((Au_{n-1},Bu_{n-1}),(Cu_{n-1},Du_{n-1}))\|~\leq~\tau~\|u_n-u_{n-1}\|.
\end{eqnarray}
Using Algorithm 4.3, and conditions (i) and (v), we get
\begin{eqnarray*}
&&\|F(v_n,w_n)-F(v_{n-1},w_{n-1})\|~\leq~\|F(v_n,w_n)-F(v_{n-1},w_{n})\|+\|F(v_{n-1},w_n)-F(v_{n-1},w_{n-1})\|\\
&&\hspace{3.9cm}~\leq~\epsilon_1\|v_n-v_{n-1}\|~+~\epsilon_2\|w_n-w_{n-1}\|\\
&&\hspace{3.9cm}~\leq~\epsilon_1\left(1+\frac{1}{n}\right)~{\cal{D}}(S(u_{n}),S(u_{n-1}))+\epsilon_2\left(1+\frac{1}{n}\right)~{\cal{D}}(T(u_{n}),T(u_{n-1}))
\end{eqnarray*}
\vspace{-0.7cm}
\begin{eqnarray}
&&\hspace{1.6cm}~\leq~\left(~\epsilon_1 l_1\left(1+\frac{1}{n}\right)~+~\epsilon_2 l_2\left(1+\frac{1}{n}\right)\right)~\|u_n-u_{n-1}\|.
\end{eqnarray}
Using conditions (iv), we get
\begin{eqnarray*}
&&\hspace{-0.7cm}\langle F(v_n,w_n)-F(v_{n-1},w_{n-1}),~J_q(H((Au_{n},Bu_{n}),(Cu_{n},Du_{n}))-H((Au_{n-1},Bu_{n-1}),(Cu_{n-1},Du_{n-1})))\rangle\\
&&\hspace{0.1cm}\leq~\langle F(v_n,w_n)-F(v_{n-1},w_{n}),~J_q(H((Au_{n},Bu_{n}),(Cu_{n},Du_{n}))-H((Au_{n-1},Bu_{n-1}),(Cu_{n-1},Du_{n-1})))\rangle\\
&&\hspace{0.5cm}+~\langle F(v_{n-1},w_n)-F(v_{n-1},w_{n-1}),~J_q(H((Au_{n},Bu_{n}),(Cu_{n},Du_{n}))-H((Au_{n-1},Bu_{n-1}),(Cu_{n-1},Du_{n-1})))\rangle\\
&&\hspace{0.1cm}\leq~(\sigma+\delta)~\|H((Au_{n},Bu_{n}),(Cu_{n},Du_{n}))-H((Au_{n-1},Bu_{n-1}),(Cu_{n-1},Du_{n-1}))\|^q
\end{eqnarray*}
\vspace{-0.7cm}
\begin{eqnarray}
&&\hspace{-9.5cm}\leq~(\sigma+\delta)\tau^q~\|u_n-u_{n-1}\|^q.
\end{eqnarray}
From (4.10)-(4.13), we have
\begin{eqnarray*}
&&\hspace{-0.6cm}\|H((Au_{n},Bu_{n}),(Cu_{n},Du_{n}))-H((Au_{n-1},Bu_{n-1}),(Cu_{n-1},Du_{n-1}))-\rho(F(v_n,w_n)-F(v_{n-1},w_{n-1})\|
\end{eqnarray*}
\vspace{-0.7cm}
\begin{eqnarray}
&&\hspace{-1.1cm}~\leq~\sqrt[q]{\tau^q+c_q \rho^q\left(\epsilon l_1\left(1+\frac{1}{n}\right)+\epsilon_2 l_2\left(1+\frac{1}{n}\right)\right)^q-\rho q (\sigma+\delta) \tau^q}~\|u_n-u_{n-1}\|.
\end{eqnarray}
Combining (4.8), (4.9) and (4.14), we have
\begin{eqnarray*}
&&\hspace{-1cm}\|u_{n+1}-u_{n}\|~\leq~\|R^{H((.,.),(.,.))}_{\rho,~M(.,.)}(z_{n+1})-R^{H((.,.),(.,.))}_{\rho,~M(.,.)}(z_{n})\|
\end{eqnarray*}
\vspace{-0.7cm}
\begin{eqnarray}
&&\hspace{1.4cm}~\leq~\theta_n~\|u_{n}-u_{n-1}\|~+~\frac{1}{r+\rho m}~\|e_{n}-e_{n-1}\|,
\end{eqnarray}
where
\begin{eqnarray}
&&\hspace{-1.1cm}~\theta_n~=~\frac{1}{r+\rho m}\sqrt[q]{\tau^q+c_q \rho^q\left(\epsilon_1 l_1\left(1+\frac{1}{n}\right)+\epsilon_2 l_2\left(1+\frac{1}{n}\right)\right)^q-\rho q (\sigma+\delta)\tau^q}.
\end{eqnarray}
Let
\begin{eqnarray}
&&\hspace{-0.5cm}~\theta~=~\frac{1}{r+\rho m}\sqrt[q]{\tau^q+c_q \rho^q(\epsilon_1 l_1+\epsilon_2 l_2)^q-\rho q (\sigma+\delta)\tau^q}.
\end{eqnarray}

Then we know that $\theta_n\to \theta$ as $n\to \infty$.
By (4.5), we know that $0 < \theta < 1$ and hence there exist $n_0 > 0$ and $\theta_0\in(0,1)$ such that $\theta_n\leq \theta_0$ for all $n\geq n_0$.
Therefore, by (4.15), we have
\begin{eqnarray}
&&\hspace{-1cm}\|u_{n+1}-u_{n}\|~\leq~\theta_0~\|u_{n}-u_{n-1}\|~+~\frac{1}{r+\rho m}\|e_{n}-e_{n-1}\|~~\forall~n\geq n_0.
\end{eqnarray}
(4.18) implies that
\begin{eqnarray}
&&\hspace{-1cm}\|u_{n+1}-u_{n}\|~\leq~\theta_0^{n-n_0}~\|u_{n_0+1}-u_{n_0}\|~+~\frac{1}{r+\rho m} \sum\limits^{n-n_0}_{j=1}\theta_0^{j-1}t_{{n}-({n-1})},
\end{eqnarray}
where $t_n=\|e_{n}-e_{n-1}\|$ for all $n\geq n_0$. Hence, for any $m\geq n> n_0$, we have
\begin{eqnarray*}
&&\hspace{-7cm}\|u_{m}-u_{n}\|~\leq~\sum\limits_{p=n}^{m-1}\|u_{p+1}-u_{p}\|
\end{eqnarray*}
\vspace{-0.7cm}
\begin{eqnarray}
&&\hspace{0.5cm}~\leq~\sum\limits_{p=n}^{m-1}\theta_0^{p-n_0}\|u_{{n_0}+1}-u_{n_0}\|+\frac{1}{r+\rho m}\sum\limits_{p=n}^{m-1}~\theta_0^{p}~\sum\limits^{p-n_0}_{j=1}\left[\frac{t_{{p}-({j-1})}}{\theta_0^{{p}-({j-1})}}\right].
\end{eqnarray}
Since $\sum\limits_{j=1}^{\infty}\|e_{j}-e_{j-1}\|~\varpi^{-j}~<~\infty,~\forall~\varpi\in~(0,1)$ and $0<\theta_0<1$, it follows that $\|u_{m}-u_{n}\|\to 0$ as $n\to \infty$, and so $\{u_n\}$ is a Cauchy sequence in $X$. From (4.6) and (4.7), it follows that $\{v_n\}$ and $\{w_n\}$ are also Cauchy sequences in $X$. Thus, there exist $u$, $v$ and $w$ such that $u_n\to u$, $v_n\to v$ and $w_n\to w$ as $n\to \infty$. In the sequel, we will prove that $v\in S(u)$. In fact, since $v_n\in S(u_n)$, we have
\begin{eqnarray*}
&&d(v,~S(u))~\leq~\|v-v_n\|~+~d(v_n,~S(u))\\
&&\hspace{1.6cm}~\leq~\|v-v_n\|~+~{\cal{D}}(S(u_n),~S(u))\\
&&\hspace{1.6cm}~\leq~\|v-v_n\|~+~\rho~\|u_n-u\|\to~0,~{\rm as}~ n\to
\infty,
\end{eqnarray*}
which implies that $d(v,~S(u))~=~0.$ Since $S(u)\in CB(X)$, it
follows that $v\in S(u)$. Similarly, it is easy to see that $w\in
T(u)$.

\vspace{0.25cm}By the continuity of $R_{\rho,~M(.,.)}^{H((.,.),(.,.))},~A,~B,~C,~D,~S,~T~{\rm and}~F$
and Algorithms 4.3, we know that $u,~v~{\rm and}~ w$ satisfy

\begin{eqnarray*}
&&\hspace{-0.6cm}u~=~R^{H((.,.),(.,.))}_{\rho,M(.,.)}~[H((Au,Bu),(Cu,Du))-\rho F(u,z)+\rho\omega].
\end{eqnarray*}
By Lemma 4.1, $(u,v,w)$ is a solution of the
problem (4.1). This completes the proof   $\square$

\vspace{.25cm}
The following example shows that assumptions (i) to (vi) of Theorem 4.6 are satisfied
for variational inclusion problem (4.1).
\begin{example}
{\rm Let $q=2$ and $X~=~\mathbb{R}^2$ with usual inner product.

\vspace{.25cm}
\noindent(i) Let $S,T:\mathbb{R}^2\multimap {\mathbb{R}^2}$ are identity mappings, then $R,S$ are $n$-Lipschitz continuous for $n=1,2$.

\vspace{.25cm}\noindent Let $A,B,C,D:\mathbb{R}^2\to \mathbb{R}^2$ be defined by
\begin{eqnarray*}
Ax=\binom{\frac{1}{10}x_1}{\frac{1}{10}x_2},~Bx=\binom{-\frac{1}{5}x_1}{-\frac{1}{5}x_2},~Cx=\binom{2x_1}{2x_2},~Dx=\binom{x_1}{x_2},~~\forall~x=(x_1,x_2)\in\mathbb{R}^2.
\end{eqnarray*}

\vspace{.25cm} \noindent
Suppose that $H:(\mathbb{R}^2\times\mathbb{R}^2)\times(\mathbb{R}^2\times\mathbb{R}^2)\to \mathbb{R}^2$ is defined by
\begin{eqnarray*}
H((Ax,By),(Cx,Dy))=Ax+Bx+Cx+Dx,~~~~~\forall~x\in\mathbb{R}^2.
\end{eqnarray*} Then, it is easy to cheek that

\vspace{.25cm}
\noindent(ii) $H((.,.),(.,.))$ is $(10,2)$-strongly mixed
cocoercive with respect to $(A,C)$ and $(5,1)$-relaxed mixed
cocoercive with respect to $(B,D)$, and $A$ is $\frac{1}{n}$-expansive for $n=10,11$ and $B$ is $\frac{1}{n}$-Lipschitz
continuous for $n=4,5$.

\vspace{.25cm} \noindent(iii) $H((A,B),(C,D))$ is $\frac{29}{n}$-mixed Lipschitz continuous with respect to $A,B,C$ and $D$ for
$n=9,10$.

\noindent Let $f,g:\mathbb{R}^2\to¨\mathbb{R}^2$ be defined by
\begin{eqnarray*}
f(x)=\binom{\frac{1}{2}x_1-\frac{4}{3}x_2}{\frac{4}{3}x_1+\frac{1}{2}x_2},~g(x)=\binom{\frac{1}{4}x_1-\frac{3}{4}x_2}{\frac{3}{4}x_1+\frac{1}{4}x_2},~\forall~x=(x_1,x_2),\in\mathbb{R}^2.
\end{eqnarray*}
\vspace{.25cm} \noindent
Suppose that $M:(\mathbb{R}^2\times\mathbb{R}^2)\to \mathbb{R}^2$ is defined by
\begin{eqnarray*}
M(fx,gx)=fx-gx,~~~~~\forall~~x=(x_1,x_2),\in\mathbb{R}^2.
\end{eqnarray*}
Then, it is easy to check that $M(f,g)$ is $\frac{1}{n}$-strongly accretive with respect to $f$ for $n=2,3$
and $\frac{1}{n}$-relaxed accretive with respect to $g$ for $n=3,4$. Moreover, for $\rho = 1$, $M$ is
generalized $\alpha\beta$-$H((.,.),(.,.))$-mixed accretive with respect to $(A,C)$, $(B,D)$ and $(f,g)$.

\vspace{.25cm}
\noindent
Let $F:\mathbb{R}^2\times \mathbb{R}^2\to \mathbb{R}^2$ are defined by
\begin{eqnarray*}
F(x,y) =~\frac{x}{4}+\frac{y}{5},~~~~~\forall~x,~y,\in\mathbb{R}^2.
\end{eqnarray*}

\noindent Then, it is easy to check that

\vspace{.2cm} \noindent(iv) $F$ is is $\frac{29}{n}$-strongly accretive with respect to $S$ and $H((A,B),(C,D))$ in the first argument for $n=30,40$ and $\frac{29}{n}$-strongly accretive with respect to $T$ and $H((A,B),(C,D))$ in the second argument for $n=40,50$;

\vspace{.25cm} \noindent(v) $F$ is is $\frac{1}{n}$-Lipschitz
continuous in the first argument for $n=3,4$ and
$~\frac{1}{n}$-Lipschitz continuous in the second argument for
$n=4,5$.

%
\vspace{.25cm}
\noindent Therefore, for the constants
\begin{eqnarray*}
&&\hspace{-.5cm}l_1=l_2=~1,~\mu_1~=~10,~\gamma_1~=~2,~\mu_2~=5,~\gamma_2~=~1,~\alpha_1~=~0.1,~\beta_1~=~0.2,\\
&&\hspace{-.5cm}\alpha~=~0.5,~\beta~=~0.25,~\sigma~=~0.725,~\delta~=~0.580,~\epsilon_1~=~0.25,~\epsilon_2~=~0.2,~\tau~=~2.9,\\
&&\hspace{-.5cm}q=2,~r~=~2.9,~m~=~0.25.
\end{eqnarray*}
obtained in (i) to (viii) above, all the conditions of the Theorem 4.7 is
satisfied for the generalized mixed variational inclusion problem
(4.1) for $\rho~=~0.35$ and $c_q=1$.}
\end{example}

\begin{remark}
{\rm If the set-valued mapping $M(f,g)$ is $\eta$-accretive and $\eta(x,y)$
is Lipschitz continuous, then $M$ becomes a new
generalized  $\alpha\beta$-$H((.,.),(.,.))$-$\eta$-mixed accretive mapping, see e.g. \cite{X-W}. We leave the proofs to
readers who are interested in this area.}
\end{remark}

\noindent{\bf Appendix: Verification (Calculations) of Example 4.7}

\vspace{.25cm}
\noindent{\bf Example 4.7}~{\rm Let $q=2$ and $X~=~\mathbb{R}^2$ with usual inner product.

\vspace{.25cm}
\noindent(i) Let $S,T:\mathbb{R}^2\multimap {\mathbb{R}^2}$ are identity mappings, then $R,S$ are $n$-Lipschitz continuous for $n=1,2$.

\vspace{.25cm}\noindent~(ii) Let $A,B,C,D:\mathbb{R}^2\to \mathbb{R}^2$ be defined by
\begin{eqnarray*}
Ax=\binom{\frac{1}{10}x_1}{\frac{1}{10}x_2},~Bx=\binom{-\frac{1}{5}x_1}{-\frac{1}{5}x_2},~Cx=\binom{2x_1}{2x_2},~Dx=\binom{x_1}{x_2},~~\forall~x=(x_1,x_2)\in\mathbb{R}^2.
\end{eqnarray*}

\vspace{.25cm} \noindent
Suppose that $H:(\mathbb{R}^2\times\mathbb{R}^2)\times(\mathbb{R}^2\times\mathbb{R}^2)\to \mathbb{R}^2$ is defined by
\begin{eqnarray*}
H((Ax,By),(Cx,Dy))=Ax+Bx+Cx+Dx,~~~~~\forall~x\in\mathbb{R}^2.
\end{eqnarray*} Then
\begin{eqnarray*}
&&\hspace{-2.5cm}\langle H((Ax,u),(Cx,u))-H((Ay,u),(Cy,u)),x-y\rangle\\
&&\hspace{0.0cm}~=~\langle (Ax+u+Cx+u)-(Ay+u+Cy+u),x-y\rangle\\
&&\hspace{0.0cm}~=~\langle Ax-Ay,x-y\rangle+\langle Cx-Cy,x-y\rangle\\
&&\hspace{0.0cm}~=~\langle (\frac{1}{10}x_1-\frac{1}{10}y_1,~\frac{1}{10}x_2-\frac{1}{10}y_2),(x_1-y_1,x_2-y_2)\rangle\\
&&\hspace{0.4cm}~+~\langle(2x_1-2y_1,2x_2-2y_2),(x_1-y_1,x_2-y_2)\rangle
\end{eqnarray*}
\vspace{-.7cm}
\begin{eqnarray}
&&\hspace{-3.7cm}~=~\frac{1}{10}\|x-y\|^2+2\|x-y\|^2
\end{eqnarray}
and
\begin{eqnarray*}
&&\hspace{0.4cm}\|Ax-Ay\|^2=~\langle Ax-Ay,Ax-Ay\rangle\\
&&\hspace{2.0cm}~=~\langle (\frac{1}{10}x_1-\frac{1}{10}y_1,~\frac{1}{10}x_2-\frac{1}{10}y_2),(\frac{1}{10}x_1-\frac{1}{10}y_1),(\frac{1}{10}x_2-\frac{1}{10}y_2)\rangle\\
&&\hspace{2.0cm}~=~\frac{1}{100}\{(x_1-y_1)^2+(x_2-y_2)^2\}
\end{eqnarray*}
\vspace{-.7cm}
\begin{eqnarray}
\hspace{-5.2cm}&&~=~\frac{1}{100}\|x-y\|^2.
\end{eqnarray}
From (4.22) and (4.23), we have
\begin{eqnarray}
&&\vspace{1.0cm}\langle H((Ax,u),(Cx,u))-H((Ay,u),(Cy,u)),x-y\rangle~\geq~10\|Ax-Ay\|^2+2\|x-y\|^2.
\end{eqnarray}
Let
\begin{eqnarray*}
&&\hspace{-2.5cm}\langle H((u,Bx),(u,Dx))-H((u,By),(u,Dy)),x-y\rangle\\
&&\hspace{0.0cm}~=~\langle (u+Bx+u+Dx)-(u+By+u+Dy),x-y\rangle\\
&&\hspace{0.0cm}~=~\langle Bx-By,x-y\rangle+\langle Dx-Dy,x-y\rangle\\
&&\hspace{0.0cm}~=~\langle (-(\frac{1}{5}x_1-\frac{1}{5}y_1),-(\frac{1}{5}x_2-\frac{1}{5}y_2)),(x_1-y_1,x_2-y_2)\rangle\\
&&\hspace{0.4cm}~+~\langle(x_1-y_1,x_2-y_2),(x_1-y_1,x_2-y_2)\rangle
\end{eqnarray*}
\vspace{-.7cm}
\begin{eqnarray}
&&\hspace{-3.7cm}~=-\frac{1}{5}\|x-y\|^2+\|x-y\|^2
\end{eqnarray}
and
\begin{eqnarray*}
&&\hspace{0.4cm}\|Bx-By\|^2=~\langle Ax-Ay,Ax-Ay\rangle\\
&&\hspace{2.0cm}~=~\langle (-(\frac{1}{10}x_1-\frac{1}{5}y_1),-(\frac{1}{5}x_2-\frac{1}{5}y_2)),(-(\frac{1}{5}x_1-\frac{1}{5}y_1),-(\frac{1}{5}x_2-\frac{1}{5}y_2))\rangle\\
&&\hspace{2.0cm}~=-\frac{1}{25}\{(x_1-y_1)^2+(x_2-y_2)^2\}
\end{eqnarray*}
\vspace{-.7cm}
\begin{eqnarray}
&&\hspace{-5.7cm}~=-\frac{1}{25}\|x-y\|^2.
\end{eqnarray}
From (4.24) and (4.25), we have
\begin{eqnarray}
&&\vspace{1.0cm}\langle H((Ax,u),(Cx,u))-H((Ay,u),(Cy,u)),x-y\rangle~\geq-5\|Ax-Ay\|^2+\|x-y\|^2.
\end{eqnarray}
From (4.22)-(4.25), $H((.,.),(.,.))$ is $(10,2)$-strongly mixed
cocoercive with respect to $(A,C)$ and $(5,1)$-relaxed mixed
cocoercive with respect to $(B,D)$, respectively.
\begin{eqnarray}
&&\hspace{-5.7cm}\|Ax-Ay\|~\geq\frac{1}{n}\|x-y\|^2,~~{\rm for}~n=10,11,\\
&&\hspace{-5.7cm}\|Bx-By\|~\leq-\frac{1}{n}\|x-y\|^2,~~{\rm for}~n=5,6.
\end{eqnarray}
From (4.27) and (4.28), $A$ is $\frac{1}{n}$-expansive for $n=10,11$ and $B$ is $\frac{1}{n}$-Lipschitz
continuous for $n=4,5$, respectively.

\noindent(iii)~Let
\begin{eqnarray*}
\|H((Ax,Bx),(Cx,Dx))-H((Ay,By),(Cy,Dy))\|^2\\
&&\hspace{-5.4cm}~=~\langle
H((Ax,Bx),(Cx,Dx))-H((Ay,By),(Cy,Dy)),\\
&&\hspace{-4.6cm}~H((Ax,Bx),(Cx,Dx))-H((Ay,By),(Cy,Dy))\rangle\\
&&\hspace{-5.4cm}~=~\langle
(Ax+Bx+Cx+Dx)-(Ay+By+Cy+Dy),\\
&&\hspace{-4.6cm}~(Ax+Bx+Cx+Dx)-(Ay+By+Cy+Dy)\rangle\\
&&\hspace{-5.4cm}~=~\langle
(\frac{29}{10}x_1-\frac{29}{10}y_1,\frac{29}{10}x_2-\frac{29}{10}y_2),(\frac{29}{10}x_1-\frac{29}{10}y_1,\frac{29}{10}x_2-\frac{29}{10}y_2)\rangle\\
&&\hspace{-5.4cm}~=~\frac{29^2}{10^2}\{(x_1-y_1)^2+(x_2-y_2)^2\},
\end{eqnarray*}
that is,
\begin{eqnarray}
\hspace{-1.5cm}\|H((Ax,Bx),(Cx,Dx))-H((Ay,By),(Cy,Dy))\|~\leq~\frac{1}{n}\|x-y\|.
\end{eqnarray}
\noindent Hence, $H((A,B),(C,D))$ is $\frac{29}{n}$-mixed Lipschitz continuous with respect to $A,B,C$ and $D$ for
$n=9,10$.\\

\noindent Let $f,g:\mathbb{R}^2\to¨\mathbb{R}^2$ be defined by
\begin{eqnarray*}
f(x)=\binom{\frac{1}{2}x_1-\frac{4}{3}x_2}{\frac{4}{3}x_1+\frac{1}{2}x_2},~g(x)=\binom{\frac{1}{4}x_1-\frac{3}{4}x_2}{\frac{3}{4}x_1+\frac{1}{4}x_2},~\forall~x=(x_1,x_2),\in\mathbb{R}^2.
\end{eqnarray*}
\vspace{.25cm} \noindent
Suppose that $M:(\mathbb{R}^2\times\mathbb{R}^2)\to \mathbb{R}^2$ is defined by
\begin{eqnarray*}
M(fx,gx)=fx-gx,~~~~~\forall~~x=(x_1,x_2),\in\mathbb{R}^2.
\end{eqnarray*}
\vspace{.25cm} Let for any $w\in \mathbb{R}^2$, we have
\begin{eqnarray*}
&&\langle M(fx,w)-M(fy,w),~x-y\rangle~=~\langle fx-w-fy+w,~x-y\rangle\\
&&\hspace{4.5cm}=~\langle fx-fy,~x-y\rangle\\
&&\hspace{4.5cm}=~\langle(\frac{1}{2}(x_1-y_1)-\frac{4}{3}(x_2-y_2)),(\frac{4}{3}(x_1-y_1)+\frac{1}{2}(x_2-y_2)),\\
&&\hspace{5.0cm}~(x_1-y_1,x_2-y_2)\rangle\\
&&\hspace{4.5cm}=~\frac{1}{2}\{(x_1-y_1)^2+5(x_2-y_2)\}^2\\
&&\hspace{4.5cm}=~\frac{1}{2}\|x-y\|^2.
\end{eqnarray*}
That is,
\begin{eqnarray}
&&\hspace{-1.0cm}\langle u-v,~x-y\rangle~\geq~\frac{1}{n}\|x-y\|^2,~\forall~x,y\in X,~u\in M(fx,w),~v\in M(fy,w).
\end{eqnarray}
Hence, $M(f,g)$ is $\frac{1}{n}$-strongly accretive with respect to $f$ for $n=2,3$ and
\begin{eqnarray*}
&&\langle M(w,gx)-M(w,gy),~x-y\rangle~=~\langle w-gx-w+gy,~x-y\rangle\\
&&\hspace{4.5cm}=-\langle gx-gy,~x-y\rangle\\
&&\hspace{4.5cm}=-\langle(\frac{1}{4}(x_1-y_1)-\frac{3}{4}(x_2-y_2)),(\frac{3}{4}(x_1-y_1)+\frac{1}{4}(x_2-y_2)),\\
&&\hspace{5.0cm}~(x_1-y_1,x_2-y_2)\rangle\\
&&\hspace{4.5cm}=-\frac{1}{4}\{(x_1-y_1)^2+(x_2-y_2)^2\}\\
&&\hspace{4.5cm}=-\frac{1}{4}\|x-y\|^2.
\end{eqnarray*}
That is,
\begin{eqnarray}
&&\hspace{-1.0cm}\langle u-v,~x-y\rangle~\geq-\frac{1}{n}\|x-y\|^2,~\forall~x,y\in X,~u\in M(w,gx),~v\in M(w,gy).
\end{eqnarray}
Hence, $M(f,g)$ is $\frac{1}{n}$-relaxed accretive with respect to $g$ for $n=3,4$.\\

\noindent From (4.31) and (4.32), $M(f,g)$ is symmetric accretive with respect to $f$ and $g$. Also for any $x\in \mathbb{R}^2$, we have
\begin{eqnarray*}
&&\hspace{-1.2cm}[H((A,B),(C,D))+\rho~M(f,g)](x)~=~[H((Ax,Bx),(Cx,Dx))+\rho~M(fx,gx)]\\
&&\hspace{3.5cm}~=~(Ax+Bx+Cx+Dx)+\rho(fx-gx)\\
&&\hspace{3.5cm}~=~(\frac{1}{10}x_1,\frac{1}{10}x_2)+(-\frac{1}{5}x_1,-\frac{1}{5}x_2)+(2x_1,2x_2)+(x_1,x_2)\\
&&\hspace{3.9cm}~+~\rho\left(\left(\frac{1}{2}x_1-\frac{4}{3}x_2,\frac{4}{3}x_1+\frac{1}{2}x_2\right)-\left(\frac{1}{4}x_1-\frac{3}{4}x_2,\frac{3}{4}x_1+\frac{1}{4}x_2\right)\right)\\
&&\hspace{3.5cm}~=~\left(\frac{29}{10}x_1,\frac{29}{10}x_2\right)+\rho\left(\frac{1}{4}x_1-\frac{7}{12}x_2,\frac{7}{12}x_1+\frac{1}{4}x_2\right)\\
&&\hspace{3.5cm}~=~\left(\left(\frac{29}{10}+\frac{1}{4}\rho\right)x_1-\frac{7}{12}\rho x_2,\left(\frac{29}{10}+\frac{7}{12}\rho\right)x_1+\frac{1}{4}\rho x_2\right),
\end{eqnarray*}
it can be easily verify that the vector on right hand side generate the whole $\mathbb{R}^2.$ Therefore, we have
\begin{eqnarray*}
&&\hspace{-1.5cm}[H((A,B),(C,D))+\rho~M(f,g)]\mathbb{R}^2~=\mathbb{R}^2.
\end{eqnarray*}
Hence, $M$ is generalized $\alpha\beta$-$H((.,.),(.,.))$-mixed accretive with respect to $(A,C)$, $(B,D)$ and $(f,g)$.

\vspace{.25cm}
\noindent (iv)
Let $F:\mathbb{R}^2\times \mathbb{R}^2\to \mathbb{R}^2$ are defined by
\begin{eqnarray*}
F(x,y) =~\frac{x}{4}+\frac{y}{5},~~~~~\forall~x,~y,\in\mathbb{R}^2.
\end{eqnarray*}
Let $u\in \mathbb{R}^2$. Since $S,T$ are Identity map, we have
\begin{eqnarray*}
&&\langle F(x,u)-F(y,u),~H((Ax,Bx),(Cx,Dx))-H((Ay,By),(Cy,Dy))\rangle\\
&&\hspace{4.5cm}=\langle \frac{1}{4}x+\frac{1}{5}u-\left(\frac{1}{4}y+\frac{1}{5}u\right),\frac{29}{10}x-\frac{29}{10}y\rangle\\
&&\hspace{4.5cm}=\langle \frac{1}{4}(x_1,x_2)-\frac{1}{5}(y_1,y_2),\frac{29}{10}(x_1,x_2)-\frac{29}{10}(y_1,y_2)\rangle\\
&&\hspace{4.5cm}=\langle \left(\frac{1}{4}x_1-\frac{1}{4}y_1,\frac{1}{5}x_2-\frac{1}{5}y_2\right), \left(\frac{29}{10}x_1-\frac{29}{10}y_1,\frac{29}{10}x_2-\frac{29}{10}y_2\right)\rangle\\
&&\hspace{4.5cm}=\frac{29}{40}\{(x_1-y_1)^2+(x_2-y_2)^2\}\\
&&\hspace{4.5cm}=\frac{29}{40}\|x-y\|^2,
\end{eqnarray*}
that is,
\begin{eqnarray}
&&\hspace{-1.0cm}\langle F(x,u)-F(y,u),~H((Ax,Bx),(Cx,Dx))-H((Ay,By),(Cy,Dy))\rangle\geq\frac{29}{n}\|x-y\|^2,
\end{eqnarray}
similarly
\begin{eqnarray}
&&\hspace{-1.0cm}\langle F(u,x)-F(u,y),~H((Ax,Bx),(Cx,Dx))-H((Ay,By),(Cy,Dy))\rangle\geq\frac{29}{n}\|x-y\|^2.
\end{eqnarray}

\vspace{.2cm}\noindent From (4.32) and (4.33), we have$F$ is is $\frac{29}{n}$-strongly accretive with respect to $S$ and $H((A,B),(C,D))$ in the first argument for $n=40,41$ and $\frac{29}{n}$-strongly accretive with respect to $T$ and $H((A,B),(C,D))$ in the second argument for $n=49,50$.

\vspace{.25cm}\noindent(v) Let
\begin{eqnarray*}
&&\hspace{-0.60cm}\|F(x,u)-F(y,u)\|^2\leq \langle F(x,u)-F(y,u),~F(x,u)-F(y,u)\rangle\\
&&\hspace{2.0cm}~=~\langle\frac{1}{4}(x_1,x_2)-\frac{1}{4}(y_1,y_2),\frac{1}{4}(x_1,x_2)-\frac{1}{4}(y_1,y_2)\rangle\\
&&\hspace{2.0cm}~=~\langle(\frac{1}{4}x_1-\frac{1}{4}y_1,\frac{1}{4}x_2-\frac{1}{4}y_2),(\frac{1}{4}x_1-\frac{1}{4}y_1,\frac{1}{4}x_2-\frac{1}{4}y_2)\rangle\\
&&\hspace{2.0cm}~=~\frac{1}{4}\{(x_1-y_1)^2+(x_2-y_2)^2\}\\
&&\hspace{2.0cm}~=~\frac{1}{16}\|x-y\|^2,
\end{eqnarray*}
that is
\begin{eqnarray}
&&\hspace{-0.60cm}\|F(x,u)-F(y,u)\|\leq \frac{1}{n}\|x-y\|,
\end{eqnarray}
similarly
\begin{eqnarray}
&&\hspace{-0.60cm}\|F(u,x)-F(u,y)\|\leq \frac{1}{n}\|x-y\|.
\end{eqnarray}

From (4.33) and (4.34), $F$ is is $\frac{1}{n}$-Lipschitz
continuous in the first argument for $n=3,4$ and
$\frac{1}{n}$-Lipschitz continuous in the second argument for
$n=4,5$.

\vspace{.25cm}
\noindent (vi)~For the constants
\begin{eqnarray*}
&&\hspace{-.5cm}l_1=l_2=~1,~\mu_1~=~10,~\gamma_1~=~2,~\mu_2~=5,~\gamma_2~=~1,~\alpha_1~=~0.1,~\beta_1~=~0.2,\\
&&\hspace{-.5cm}\alpha~=~0.5,~\beta~=~0.25,~\sigma~=~0.725,~\delta~=~0.580,~\epsilon_1~=~0.25,~\epsilon_2~=~0.2,~\tau~=~2.9,\\
&&\hspace{-.5cm}q=2,~r~=~2.9,~m~=~0.25.
\end{eqnarray*}
obtained from above conditions (i) to (vi), one can easily verify, for $c_q=1$ and $\rho = 3.8$ and , that the
condition (vi) of Theorem 4.6, given by
$$0<\sqrt[q]{\tau^q+c_q \rho^q(\epsilon_1 l_1+\epsilon_2 l_2)^q-\rho q (\sigma+\delta)\tau^q}<{r+\rho m},$$
~where~$r~=~(\mu_1\alpha_1^q~-~\mu_2\beta_1^q)+(\gamma_1+\gamma_2)$ and $m~=~\alpha-\beta$,
and $\alpha>\beta$, $\mu_1>\mu_2,~\alpha_1>\beta_1$ and $\gamma_1,~\gamma_2,\rho>0~$ holds.}

\vspace{0.25cm}\noindent{\bf Acknowledgement:} This research is supported by National Board for Higher Mathematics, Mumbai, Maharashtra India (Grant number: 2/40(13)/2014/R$\&$D-II/13106).


\bibliographystyle{amsplain}

\end{document}